\documentclass[journal,twocolumn]{IEEEtran}
\usepackage{cite}
\usepackage{amsmath,amsthm,amssymb,amsfonts}
\usepackage{algorithmic}
\usepackage{graphicx}
\usepackage{textcomp}
\usepackage{xcolor}
\usepackage{thmtools}
\usepackage{mathtools}
\usepackage{bm}
\usepackage{mdframed}
\usepackage{comment}

\newtheorem{theorem}{Theorem}
\newtheorem{lemma}{Lemma}
\newtheorem{corollary}{Corollary}

\theoremstyle{definition}
\newtheorem{definition}{Definition}

\theoremstyle{remark}
\newtheorem*{remark}{Remark}

\newcommand{\dtil}[1]{\dot{\tilde{#1}}(t)}
\newcommand{\til}[1]{\tilde{#1}(t)}
\newcommand{\calE}{\mathcal{E}}
\newcommand{\calC}{\mathcal{C}}
\newcommand{\ev}[1]{\mathcal{E}_{x(#1)}}
\newcommand{\Rd}{\mathbb{R}^d}
\newcommand{\knl}{\mathfrak{K}}
\newcommand{\Hx}{\mathcal{H}_X}
\newcommand{\Ho}{\mathcal{H}_{\Omega}}
\newcommand{\Vo}{\mathcal{V}_\Omega}
\newcommand{\Hn}{\mathcal{H}_{\Omega_n}}
\newcommand{\dotp}[1]{\left\langle#1\right\rangle_{\mathcal{H}_X}}
\newcommand{\restr}{\mathbf{R}_\Omega}
\newcommand{\exten}{\mathbf{E}_X}
\newcommand{\prj}{\mathbf{P}_\Omega}
\newcommand{\calR}{\mathcal{R}}
\newcommand{\RR}{\mathbb{R}}

\begin{document}

\title{Partial Persistence of Excitation in RKHS Embedded Adaptive Estimation}

\author{Jia~Guo,~
        Sai~Tej~Paruchuri,~
        and~Andrew~J.~Kurdila
\thanks{The authors are with the Department
of Mechanical Engineering, Virginia Tech, Blacksburg,
VA, 24060 USA. 
E-mail: {\tt\small \{jguo18, saitejp, kurdila\}@vt.edu}.}
}


\maketitle

\begin{abstract}
In this paper, an adaptive non-parametric method is proposed to estimate the scalar-valued nonlinear function that appears in  uncertain systems governed by ordinary differential equations (ODEs)
. By employing an infinite-dimensional reproducing kernel Hilbert space (RKHS) as the hypothesis space, the nonlinear estimation problem in finite-dimensional Euclidean space is recast into that of constructing a  linear observer  in the  infinite-dimensional RKHS. 
The analysis of convergence is facilitated by the introduction of a novel  condition of partial persistent excitation (partial PE), 
which is defined for a  subspace $\Ho\subseteq \Hx$ of the RKHS $\Hx$. Using this condition, we prove that the projection of the function estimation error onto the PE subspace $\Ho$ converges in norm asymptotically to zero. 
While this is an abstract notion of convergence that depends implicitly on the kernel used to define the RKHS, we derive conditions that ensure the pointwise convergence of the function estimates over the PE subset $\Omega$. This paper additionally introduces a weaker but geometrically intuitive notion of a partial PE condition, one that  resembles PE conditions as they have been formulated historically in Euclidean spaces. Sufficient conditions are derived that describe when the two conditions are equivalent. Finally, qualitative properties of the convergence proofs derived in the paper are illustrated with numerical simulations.

\end{abstract}



%
%
\section{Introduction}
\label{sec:intro} 
Adaptive estimation is now a common approach for online identification of uncertain systems governed by ordinary differential equations (ODEs). In this paper we study a model problem that is governed by the equation 
\begin{align}
\label{eq:new1}
\dot{x}(t)&=Ax(t)+B f(x(t))
\end{align}
where the state $x(t)\in X= \Rd$, the matrix $A\in \RR^{d\times d}$ and the vector $B\in \RR^{d\times 1}$ are known, and the unknown is the typically nonlinear function $f:\RR^d \rightarrow \RR$. 

In this paper, the unknown function $f$ is assumed to be an element contained in $\Hx$, a generally infinite-dimensional reproducing kernel Hilbert space (RKHS) of real-valued functions over $\Rd$. We define an adaptive estimation scheme that generates state estimates $\hat{x}(t)\in \RR^d$ and function estimates $\hat{f}(t):=\hat{f}(t,\cdot)\in \Hx$ of the unknown function $f\in \Hx$ at each time $t\in \RR^+$ using the state history $\{x(\tau)\}_{\tau\leq t}$. The estimate $\hat{f}:\RR^+\times \Rd \rightarrow \RR$ is a function of both time and space. Here and in the remainder of this paper we use the notation $\hat{f}(t):=\hat{f}(t,\cdot)$ as a shorthand notation for the spatial function $x\mapsto \hat{f}(t,x)$ at fixed time $t\in \RR^+$.   
If we denote the state error by $\tilde{x}(t)=x(t)-\hat{x}(t)$ and the function estimate error by $\tilde{f}(t)=f-\hat{f}(t,\cdot)$, the RKHS embedding method generates errors that satisfy the equations
\begin{equation}
\label{eq:err_rkhs0a}
    \begin{bmatrix}\dtil{x} \\ \dtil{f} \end{bmatrix} 
    = \begin{bmatrix} A & B\ev{t} \\ -\mu(B\ev{t})^*P & 0 \end{bmatrix} 
    \begin{bmatrix}\til{x} \\ \til{f} \end{bmatrix}
\end{equation}
which resemble the error equations in the conventional adaptive estimator. Here $\calE_{x(t)}$ is the evaluation functional at $x(t)$, and $\calE^*_{x(t)}$ is its adjoint. The evaluation functional is defined  via the identity  $\calE_{x}f=f(x)$ for any $x\in \RR^d$ and $f\in \Hx$. Constant $\mu>0$ is the learning rate, and matrix $P$ satisfies a Lyapunov equation as summarized in Section \ref{sec:problem}. 
In contrast to conventional approaches to adaptive estimation of unknown functions in uncertain ODEs, these error equations define a distributed parameter system (DPS) with trajectories that evolve in $\RR^d\times \Hx$. In some sense, these equations can be viewed as a description of the limiting behavior of adaptive estimation methods in $\Rd\times\RR^n$ (as in Equation \ref{eq:eqconv}), when the number $n$ of the parameters used to describe the function estimate   increases to infinity. 

The primary contribution of this paper is the derivation of conditions under which the function error $\tilde{f}(t)$ converges to zero in $\Hx$-norm as $t\rightarrow \infty$. Our convergence analysis depends in a critical way on the appropriate definition  of a  closed subspace $\Ho \subseteq \Hx$  that is partially persistently excited as an RKHS subspace of $\Hx$.  
With these definitions we obtain an entire family of convergent estimation methods that are not in general equivalent: they depend on the topology of the RKHS $\Hx$.  We show that when $\Ho$ is a partially persistently excited subspace of $\Hx$, the RKHS norm $\|\cdot\|_{\Hx}$ of the projection $\mathbf{P}_\Omega\tilde{f}(t,\cdot)$ converges to zero as $t\rightarrow \infty$. In this way we show that the partial PE condition in the RKHS embedding method results in a guarantee of convergence that is qualitatively similar to the guarantee of parametric convergence that the conventional partial PE condition ensures for adaptive estimation in Euclidean spaces.

%
\subsection{Historical Perspectives}

To motivate our approach, we review some well-known techniques for estimating uncertain ODE systems that have the form of Equation \eqref{eq:new1}. Over the years, the most well-understood class of such systems are those in which the unknown function $f$ is a linear function of finitely many real parameters. The associated estimation problem is then one of (real) parameter estimation. For general nonlinear systems, it is a common practice to assume the unknown function $f$ appearing in Equation \eqref{eq:new1} can be written as a superposition of a finite set of nonlinear regressor functions $\{\phi_i(x)\}_{i=1}^n$, which is often called the linear-in-parameter (LIP) assumption \cite{PoFar}. 
In particular, it is fairly standard to assume that $f$ can be approximated by the LIP representation with small uniform error over a compact subset in which the unknown system trajectory $x(t)$ evolves \cite{PoFar}. To reduce the amount of prior information that is needed  regarding the form of the unknown function, appeals to approximation theory have resulted in techniques that are  systematic in  selecting regressor functions. For example, universal regressors such as radial basis functions (RBFs) have been incorporated in online adaptive estimation and control since 1990s \cite{slotine1991,slotine1992}. A more comprehensive introduction on this topic is given in \cite{PoFar}. 

Once the regressors $\Phi=\{\phi_{i}\}_{i=1}^n$ are chosen, the function $f$ is assumed to lie in the span of the regressors $\mathcal{H}_n:=\text{span}\{\phi_i : 1\leq i\leq n\}$. Then the adaptive learning law is introduced for estimating  $\hat{\alpha}(t)=\{\hat{\alpha}_i(t)\}_{i=1}^n$ {\em of the parameters} $\alpha=\{\alpha_i\}_{i=1}^n$ of the LIP representation, and this law can be selected from a host of alternatives described in classical texts like \cite{sastry2011book,naranna,IaSu,PoFar}.  One example of the error equations for this learning law can be written in the form 
\begin{equation}
    \begin{bmatrix}\dtil{x} \\ \dtil{\alpha} \end{bmatrix} 
    = \begin{bmatrix} A & B\Phi^T(x(t)) \\ -\mu\Phi(x(t))B^TP & 0 \end{bmatrix} 
    \begin{bmatrix}\til{x} \\ \til{\alpha} \end{bmatrix} \label{eq:eqconv}
\end{equation}
where $A,B,\tilde{x}(t),\mu,$ and $P$ are defined as in Equation \eqref{eq:err_rkhs0a} and $\tilde{\alpha}(t)=\alpha-\hat{\alpha}(t)$ is the error in the parameter estimates. These equations define an evolution law in the Euclidean space $\RR^d\times \RR^n$. 
The convergence of this formulation, \textit{in time} for a fixed number of parameters $n$, is originally proven using a PE condition in classical references such as \cite{narendra1977b,anderson1977exponential}. 

Although the stability and convergence of similar online adaptive estimation formulations has been studied extensively \cite{panteley2003,annaswamy2018,sastry2019,annaswamy2019}, it is recognized that the fidelity of the estimate $\hat{f}(t,\cdot)=\sum_{i=1}^n\phi_i(\cdot) \hat{\alpha}_i(t)$ of the unknown function $f$ is often important in its own right. In other words, it is of interest to understand how well $\hat{f}(t,\cdot)$ approximates $f$.  Most of the discussions about adaptive systems, in the sense described in classical texts like \cite{sastry2011book,naranna,IaSu,PoFar},  focus on the convergence of the  parameter estimates in $\RR^n$, while the convergence of the associated function estimates {\em in space} is rarely discussed in terms of norms on function spaces. Often, such convergence in the function estimate is described only in qualitative terms.  As an example, a rough assessment of the reliable region over which function estimates are made is given in \cite{wang2006} and the following work in \cite{wu2019}. Likewise, a qualitative illustration of the convergence of error over trajectories is given in the example in \cite{farrell99wavelet}. But again, there is seldom result on the rigorous connection between the range of the  trajectory $t\mapsto x(t)$  and the convergence of function estimates {\em in some norm topology on functions} for ODE systems. It would be a significant improvement over existing methods formulated in Euclidean space to generate concrete bounds on the asymptotic behavior of the estimates in some norm $\|\cdot\|_{\Hx}$ on functions. One of the goals of this paper is to establish such convergence. 


%
\subsection{Non-parametric Estimation of Functions}

To tease out a description of the function estimate and its properties, it can be useful to consider the  finite-dimensional hypothesis space $\mathcal{H}_n = span\{\phi_i\}_{i=1}^n$  as an approximation of some infinite-dimensional, general function space $\Hx$ to which the unknown function belongs. The important question then becomes how to choose the appropriate infinite dimensional space. Recent progress in both approximation and learning theory focuses on the role of RKHS as a means of solving certain approximation problems that take as data scattered samples  \cite{wendland,rasmussen}. As further explained in Section \ref{sec:RKHS_review}, RKHS have several attractive  properties that are suitable for extending the classical online adaptive methodologies to the infinite dimensional setting. Many different types of function spaces can be defined or approximated  by  selection of a  suitable  reproducing kernel.
It is also worth-noting that any RKHS can be constructed as the completion of the space spanned by kernel basis functions whose centers are located in the domain \cite{aronszajn1950,berlinet}. This fact leads to convenient and pragmatic schemes to parameterize the abstract form of function estimates and analyze the approximation error. As we discuss more fully below, it is straightforward to define the finite dimensional spaces of approximants $\mathcal{H}_n$ directly in terms of the translates of the kernel that defines the infinite dimensional space $\Hx$. 

The problem of identifying or approximating an unknown function directly from observations utilizing an RKHS setting  has been studied in many references. 
A good and popular account of the theory from an approximation theory viewpoint is found in \cite{wendland}. Another well-cited account, one with a foundation in machine learning strategies, is \cite{rasmussen}. This latter reference is particularly useful in its description of RKHS in problems of Bayesian estimation, see for example \cite{chowdhary2014}. The RKHS framework has come to widespread use in the general area of estimation of Gaussian processes.  For example, a  probabilistic method based on Gaussian process regression is proposed in \cite{Pillonetto2010a,Pillonetto2011b}, where the impulse response of an unknown SISO system is assumed to be a realization of some Gaussian process. The \textit{a posteriori} process is computed from the observed data using Gaussian process regression. Then finite-dimensional approximations are obtained by a projection that minimizes the mean square error. Using this approach, the selection and design of kernel (covariance) functions are studied in order to reduce the approximation error in different situations \cite{Pillonetto2011, Chen2018, Libera2019}.


This paper concentrates on how adaptive estimation can be formulated in a way to exploit the advantages of an RKHS setting. It builds on some recent work for online adaptive estimation of nonlinear ODEs in continuous time and seeks to address some unanswered questions suggested there. 
The framework for RKHS embedded adaptive estimation has been  proposed in \cite{bmpkf2017C,bmpkf2017} as an initial attempt.
This  method is inspired by the earlier work in \cite{bsdr1997,bdrr1998,d1993,dr1994,dr1994pe}, which studies the adaptive estimation problem for certain types of DPS defined using  a pivot space structure in Banach and Hilbert spaces. Some key, but initial,  theoretical questions regarding the RKHS embedding method for uncertain ODEs have already been studied and documented in \cite{bmpkf2017C,bmpkf2017}.  
Reference \cite{kl2013} is an early example that applies the RKHS embedding strategy to $\mathcal{L}^1$-adaptive control problems. The basic method is further extended to solve consensus estimation problems \cite{bobade2}.

\subsection{Our Contributions}

In this paper two definitions are introduced of what constitute a partially PE subspace $\Ho\subset\Hx$ and the corresponding PE subset $\Omega\subset X$. It should be noted that these two definitions of partial PE can be viewed as the infinite dimensional analogs of two classical notions of PE in Euclidean spaces studied in \cite{narendra1977b}. The central theorem proven in this paper states a general result that, when the new partial PE conditions holds over a subset $\Omega\subset X$, we have $\|\prj(f-\hat{f}(t))\|_{\Hx}\rightarrow 0$ as $t\rightarrow \infty$ where $\prj:\Hx \rightarrow \Ho$ is the orthogonal projection onto the partially PE subspace $\Ho$. This result enables us to conclude the local convergence of the RKHS embedding method: we have $|\mathbf{P}_{\Omega}(f-\hat{f}(t))(x)|\rightarrow 0$ for each $x\in \Omega$ as $t\rightarrow \infty$, and this convergence is uniform in $x\in \Omega$ if the kernel is uniformly bounded. This precise characterization of local convergence for adaptive estimation is the main contributions of this paper. An additional contribution of this paper.  The newly proposed PE condition (PE.1) has a different structure from  those in the our previous studies \cite{kur95,kgp2019}, which take the form of PE.2 as in Definition \ref{def:PE2}. The version of the condition in (PE.2) is not always equivalent to that in (PE.1). In this sense, the situation resembles that in Euclidean space described in \cite{narendra1977b}. The equivalence between the two formulations is proven in a new fashion that is different than that in the classical analysis in \cite{narendra1977b}. We establish equivalence under the condition that the functions in the unit ball of $\Hx$ are uniformly equicontinuous. A criterion for selecting kernels is described in the numerical examples, and  this strategy  bridges the constructive conclusions about PE in \cite{kur95,kgp2019,pgk2020a} and the local convergence proven in this paper. 

The organization of this paper is as follows. The theory of RKHS is  briefly reviewed in Section \ref{sec:RKHS_review}. The review is mainly focused on the case that the RKHS is  uniformly embedded in the space of continuous functions. The structure and properties of subspaces in the RKHS are discussed in detail. The problem setup is introduced in Section \ref{sec:problem}, including the statement of the governing equations of the RKHS embedded adaptive estimation and the new definitions of partial PE conditions. The main results guaranteeing  local convergence are stated in Section \ref{sec:convergence}. The exploration of  the implications of the new partial PE conditions is contained in Section \ref{sec:equiv_PE}. Numerical simulations are presented  in Section \ref{sec:num_sim}.

\section{Reviews of RKHS}
\label{sec:RKHS_review}
\subsection{Notations}
In this paper, we use $X:=\Rd$ to denote the state space. Norms on different spaces are distinguished using subscripts: $\|\cdot\|$ for the Euclidean spaces $\Rd$ (state space) or $\mathbb{R}^n$ (parameter space), $\|\cdot\|_{op}$ for the uniform operator norm, and $\|\cdot\|_{\Hx}$ for the RKHS $\Hx$. The inner product on $\Rd$ and $\Hx$ are written as $(\cdot,\cdot)$ and $\dotp{\cdot,\cdot}$ respectively.

The subscription $\Omega$ in $\Ho$ denotes the index set of the RKH subspace, where $\Omega\subseteq X$ is a non-empty subset. See Corollary \ref{cor:projection} for more details. In particular, the subspace $\Ho=\Hx$ when $\Omega=X$. The evaluation functional $\calE_x:\Hx\rightarrow\mathbb{R}$ denotes the operation that evaluate a function $f$ at a point $x$ in the domain, so $\calE_x f=f(x)$. The asterisk superscript denotes the adjoint of bounded linear operators acting between  Hilbert spaces, e.g., as in  $\calE_x^*$ and $B^*$.

\subsection{RKHS Uniformly Embedded in $\mathcal{C}(X)$}
A real RKHS $\Hx$ is a Hilbert space of real-valued functions defined over $X$, which admits a reproducing kernel $\knl(\cdot,\cdot)$. The reproducing kernel $\knl(\cdot,\cdot)$ is a real-valued function defined on $X\times X$ that satisfies the following reproducing property: for any function $f\in\Hx$ and $x\in X$, the point evaluation $f(x)$ equals the inner product with the kernel basis function located at $x$, i.e.,
\begin{equation*}
    \calE_x f = f(x) = \dotp{f,\knl(\cdot,x)},
\end{equation*}
where $\calE_x$ denotes the point evaluation functional at $x\in X$. It is straightforward to show that the adjoint operator $\calE_x^*:\RR\rightarrow\Hx$ can be expressed as the multiplication operator $\calE_x^*\alpha = \alpha\knl(\cdot,x)$ for all $\alpha\in\RR$ and $x\in X$. This explicit expression is used in the implementation of the estimator equations.

All the functions in $\Hx$ are defined on the same domain $X$. For simplicity, we denote the kernel basis function $\knl(\cdot,x)$ by $\knl_x(\cdot)$, which clearly is an element in $\Hx$. Since $\Hx$ is a Hilbert space, by the Cauchy-Schwartz inequality we have for every $f\in\Hx$, 
\begin{equation}
\label{eq:reproduce}
    |\calE_xf| = \left|\dotp{f,\knl_x}\right| \leq \|f\|_{\Hx} \|\knl_x\|_{\Hx}.
\end{equation}
Since the kernel basis function $\knl_x$ belongs to $\Hx$, by the reproducing property we have 
$
    \|\knl_x\|_{\Hx} = \sqrt{\knl(x,x)}.
$
When we take the supremum for both sides in Equation \eqref{eq:reproduce} over all $x\in X$, the left hand side then becomes $\sup_{x\in X}|f(x)|$, which equals the uniform norm $\|f\|_{\infty}$. Thus we have
\begin{equation*}
    \|f\|_{\infty} = \sup_{x\in X}|f(x)|\leq \|f\|_{\Hx}\sup_{x\in X}\sqrt{\knl(x,x)}.
\end{equation*}
In the following discussions, we assume there exists a uniform bound $\bar{k}<\infty$ such that $\sqrt{\knl(x,x)}\leq \bar{k}$ for all $x\in X$. For example, this assumption holds naturally if the kernel is a radial basis function such as the Gaussian kernel, where $\knl(x,y):=\mathcal{R}(\|x-y\|)$ and $\mathcal{R}(\cdot)$ is radial. Under this assumption, for all $f\in\Hx$ we have $\|f\|_\infty\leq \bar{k}\|f\|_{\Hx}$. Then it follows that the RKHS $\Hx$ is uniformly embedded in the space of continuous functions $\mathcal{C}(X)$, that is,
\begin{equation*}
    \Hx\hookrightarrow\mathcal{C}(X).
\end{equation*}

\subsection{Subspaces of an RKHS}
It is proven in the pioneering work by Aronszajn \cite{aronszajn1950} that the span of all kernel functions $\knl_x$ for $x\in X$ forms a dense subspace in $\Hx$, in other words,
\begin{equation*}
    \Hx = \overline{span\{\knl_x:x\in X\}}.
\end{equation*}
Let $\Omega\subseteq X$ be any subset. Denote by $\restr$ the operator that restricts the domain of a function to $\Omega$,
\begin{equation*}
    \restr f := f|_\Omega.
\end{equation*}
The following theorem in \cite{aronszajn1950,berlinet} summarizes some key properties of the space $\restr\Hx$ of restricted functions. 
\begin{theorem}[Theorem 6 in \cite{berlinet}]
\label{th:ber6}
Let $\Hx$ be a Hilbert space of functions defined on $X$ with reproducing kernel $\knl$. Let $\Omega$ be a non-empty subset of $X$. The restriction  $\knl_{\Omega}:=\knl|_{\Omega\times\Omega}$ is the reproducing kernel of the space $\restr\Hx$. For $\bar{f}\in\restr\Hx$, its norm is defined as
\begin{equation}
\label{eq:rh_norm}
    \|\bar{f}\|_{\restr\Hx} := \min\{\|f\|_{\Hx}:\restr f = \bar{f},f\in\Hx\}.
\end{equation}
\end{theorem}
\noindent We define the extension operator $\exten:\restr\Hx \rightarrow \Hx$ that maps the restricted function $\bar{f}$ to the corresponding minimizer $f$ as in Equation \eqref{eq:rh_norm},
\begin{equation*}
    \exten\bar{f}:= \arg\min\{\|f\|_{\Hx}:\restr f = \bar{f},f\in\Hx\}.
\end{equation*}
It turns out that the space of extensions $\exten\restr\Hx$ form a closed subspace of $\Hx$, the properties of which are characterized in the following corollary.
\begin{corollary}
\label{cor:projection}
Suppose $\knl:X\times X\rightarrow\mathbb{R}$ is a kernel that induces an RKHS $\Hx$. Let $\Omega\subseteq X$ be a non-empty subset of the domain $X$. Then the space $\exten\restr\Hx$ is a closed subspace of $\Hx$. It is orthogonal to the closed subspace $\Vo$ defined as follows,
\begin{equation*}
    \Vo=\{v\in\Hx: v(x)=0 \; \forall x\in\Omega \}.
\end{equation*}
The total RKHS $\Hx$ can be written as the direct sum $\Hx = \exten\restr\Hx \oplus \Vo$, and the composition $\prj=\exten\restr$ is a projection operator. Moreover, we have the explicit expression of $\exten\restr\Hx$ as follows
\begin{equation*}
    \exten\restr\Hx = \Ho = \overline{span\{\knl_x:x\in\Omega\}}.
\end{equation*}
\end{corollary}
\proof{See the Appendix.}

\begin{corollary}
\label{cor:proj_equiv}
If two functions in the RKHS $f,g\in\Hx$ have same projection over a non-empty subset $\Omega\subseteq X$, then for any $x\in\Omega$, $f(x)=g(x)$.
\end{corollary}
\begin{proof}
The identity $\prj(f-g)=0$ implies that $f-g\perp\Ho$, which means $f-g\in\Vo$. Thus by the definition of $\Vo$, we have $f(x)=g(x)$ for all $x\in \Omega$.
\end{proof}

As shown in the next section, the notions of projection and subspace studied in Corollary \ref{cor:projection} come in handy when defining the partial PE conditions. The conclusion in Corollary \ref{cor:proj_equiv} about two functions having same projections will be used for depicting the local convergence of functions estimates. In reference \cite{kgp2019}, the connection between the subspaces $\Ho$ and the index set $\Omega$ is further explored, and a necessary condition of partial PE is proposed in a constructive manner.

When the index set $\Omega_n=\{x_i\}_{i=1}^N$ is finite, the subspace $\Ho=span\{\knl_{x_i}:i=1,...,N\}$ degenerates to a finite-dimensional function space where the function norm is equivalent to the Euclidean norm of coefficient vectors. In such case, the RKHS embedded estimator degenerates to the common adaptive estimator of parameters. 

%
%
\section{Problem Statement}
\label{sec:problem}
\subsection{RKHS Embedding Method}
Overall, we want to structure the notion of partial PE in the RKHS embedding method in a way that the associated guarantees of $\Hx$-norm convergence closely resemble the well-known conventional case for real-parameter estimation. For this reason, we have elected, for purposes of comparison, to summarize the conventional case for estimation in $\Rd\times\RR^n$ and subsequently carefully define the situation for RKHS embedding in $\Rd\times\Hx$. As discussed in the introduction, the unknown system is governed by a set of nonlinear ODE having the following form,
\begin{equation}
\label{eq:orig_dyn_clas}
\dot{x}(t) = Ax(t) + Bf(x(t)),
\end{equation}
where $A\in\mathbb{R}^{d\times d}$ is a known Hurwitz matrix, $B\in \mathbb{R}^{d\times 1}$ is known, and $f:\Rd\rightarrow \mathbb{R}$ is the unknown function. 

Under the LIP assumption, the unknown function has the form of $f=\sum_{i=1}^n \alpha_i^* \phi_i(\cdot)$ where $\phi_i:\mathbb{R}^d\rightarrow\mathbb{R}$ for $i=1,...,n$ denote the regressors. A classical example of the adaptive state observer and update law are 
\begin{equation}
\begin{aligned}
\label{eq:clas_estimator}
\dot{\hat{x}}(t) &= A\hat{x}(t) + B\Phi^T(x(t))\hat{\alpha}(t), \\
\dot{\hat{\alpha}}(t) &= \mu \big[B\Phi^T(x(t))\big]^TP(x(t) - \hat{x}(t)),
\end{aligned}
\end{equation}
which induces the following error equations
\begin{equation}
    \label{eq:err_rd}
    \begin{bmatrix}\dtil{x} \\ \dtil{\alpha} \end{bmatrix} 
    = \begin{bmatrix} A & B\Phi^T(x(t)) \\ -\mu\Phi(x(t))B^TP & 0 \end{bmatrix} 
    \begin{bmatrix}\til{x} \\ \til{\alpha} \end{bmatrix}.
\end{equation}
Here $\til{x}=x(t)-\hat{x}(t)\in \Rd$ is the error in state estimates $\hat{x}(t)$ of the actual states $x(t)$, and $\tilde{\alpha}(t)=\alpha^*-\hat{\alpha}(t)\in \mathbb{R}^n$ is the error in the parameter estimates $\hat{\alpha}(t)$ of the actual parameters $\alpha^*$. The matrix $P\in \mathbb{R}^{d\times d}$ is the unique positive definite solution of the Lyapunov equation $A^TP+PA=-Q$ for some fixed symmetric positive definite matrix $Q\in \mathbb{R}^{d\times d}$. The function $\Phi(\cdot)=[\phi_1(\cdot),...,\phi_n(\cdot)]^T$ is the vector of regressor functions.

The RKHS embedded adaptive estimation proposed in \cite{bmpkf2017C} relaxes the LIP assumption by assuming the unknown function $f$ belongs to an RKHS $\Hx$. The governing equations are 
\begin{equation}
\label{eq:rkhs_estimator}
\begin{aligned}
\dot{\hat{x}}(t) &= A\hat{x}(t) + B\calE_{x(t)}f, \\
\dot{\hat{f}}(t) &= \mu (B\calE_{x(t)})^*P(x(t) - \hat{x}(t)),
\end{aligned}
\end{equation}
which induce the infinite-dimensional error equations 
\begin{equation}
\label{eq:err_rkhs0}
    \begin{bmatrix}\dtil{x} \\ \dtil{f} \end{bmatrix} 
    = \begin{bmatrix} A & B\ev{t} \\ -\mu(B\ev{t})^*P & 0 \end{bmatrix} 
    \begin{bmatrix}\til{x} \\ \til{f} \end{bmatrix}.
\end{equation}
Most of the terms in Equation \eqref{eq:err_rkhs0} have the same meaning as in Equation \eqref{eq:err_rd}. The essential difference is that the ``state'' $\tilde{f}(t) = f-\hat{f}(t)$ is now a function that evolves in the RKHS $\Hx$. See \cite{bmpkf2017} for detailed discussions about the well-posedness of this equation, as well as the functional properties. If the system governed by the error equations has an asymptotically stable equilibrium at the origin in $\Rd\times\Hx$, the estimator converges in the norm on $\Rd\times\Hx$. Due to the infinite dimensionality of $\Hx$, the convergence can be established or analyzed in multiple topologies, each of which defines a notion of convergence that might be distinct. 



\subsection{Partial PE Conditions}
The usual notion of PE in Definition \ref{def:PE_rn} is defined to study the asymptotic stability of the classical error equations in Equation \eqref{eq:err_rd}. Assuming that the regressors $\Phi(\cdot)$ are piecewise continuous and Lipschitz, this condition is a sufficient and necessary condition for the uniform asymptotic stability of the error Equations \eqref{eq:err_rd}. Details are discussed in references \cite{narendra1977b,anderson1977exponential}.
\begin{definition}[PE in $\mathbb{R}^n$]
\label{def:PE_rn}
A trajectory $t\mapsto x(t)\in \Rd$ persistently excites a family of regressor functions $\Phi:\Rd\rightarrow \mathbb{R}^n$ if there exist constants $T_0,\Delta,\gamma>0$ such that
\begin{equation}
\label{eq:PE_rn}
\int_t^{t+\Delta} v^T\Phi(x(\tau))\Phi^T(x(\tau))v d\tau \geq \gamma\|v\|^2,
\end{equation}
for each $t\geq T_0$ and $v\in \mathbb{R}^n$.
\end{definition}
\noindent Sometimes the condition above is too strong to guarantee, and in this case it can be advantageous to introduce partial PE conditions. Suppose the matrix $P_m\in\mathbb{R}^{m\times n}$ where $m<n$ has full rank. Then it is a surjective linear mapping from $\mathbb{R}^n$ to $\mathbb{R}^m$. One notion of partial PE can be defined by setting $\alpha=P_mv$. In this way, the regressors in $P_m\Phi$ instead of $\Phi$ are persistently excited. Convergence of parameters can still be achieved with regard to $P_m\Phi$. Details are discussed in \cite{Boyd1983On}.

In this paper, we introduce the two definitions of the partial PE condition for the RKHS embedding method. The operator $\prj:\Hx\rightarrow \Ho$ is the orthogonal projection that maps any function $g\in \Hx$ onto the closed RKH subspace $\Ho$, which is defined in terms of the indexing set $\Omega\subseteq X$. 
\begin{definition}[PE.1 in an RKHS]
\label{def:PE1}
A trajectory $t\mapsto x(t)\in X$ persistently excites the subspace $\Ho$ of an RKHS $\Hx$ if there exist a constant $\gamma_1>0$ and time constants $T_0, \delta,$ and $\Delta$, such that for each $t\geq T_0$ and any $g\in \Hx$, there exists $s\in [t,t+\Delta]$ such that 
\begin{equation}
\label{eq:PE1}
    \left| \int_s^{s+\delta} \ev{\tau} g d\tau \right| \geq \gamma_1\|\prj g\|_{\Hx}.
\end{equation}
\end{definition}
\begin{definition}[PE.2 in an RKHS]
\label{def:PE2}
A trajectory $t\mapsto x(t)\in X$ persistently excites the subspace $\Ho$ of the RKHS $\Hx$ if there exist a constant $\gamma_2>0$ and time constants $T_0$ and $\Delta$, such that for all $t\geq T_0$ and any $g\in \Hx$,
\begin{equation}
\label{eq:PE2}
    \int_{t}^{t+\Delta}\left\langle\ev{\tau}^*\ev{\tau}g,g \right\rangle_{\Hx} d\tau \geq \gamma_2\|\prj g\|^2_{\Hx}.
\end{equation}
\end{definition}
\noindent In either of these two definitions, we say that the set $\Omega$ is partially PE when the subspace $\Ho$ is PE. These two definitions are analogous to those studied in the classical scenario in \cite{narendra1977b}, but here they are expressed in terms of the evaluation operator $\mathcal{E}_x$ rather than the regressor functions $\Phi$. The PE.1 condition is used in Theorem \ref{th:pe_uas} to prove the convergence of the estimator, while the PE.2 condition seems more intuitive and leads to constructive results \cite{kgp2019,pgk2020a}. The connection between the two conditions is explored below in Section \ref{sec:equiv_PE}. Also, Definition \ref{def:PE1} should be compared to that given in \cite{bsdr1997} and the references therein for evolution equations cast in terms of Gelfand triples. 

Carefully note that to guarantee the subspace $\Ho$ is PE in the sense of the Definitions \ref{def:PE1} and \ref{def:PE2}, all the functions $g\in \Hx$ should be tested. We denote this by saying that the pair $(\Ho,\Hx)$ is PE, where the first argument is the PE subspace and the second argument is the space in which the functions are tested. 
We will see in the Theorem \ref{th:pe_uas} that the strongest conclusion regarding convergence is achieved when the pair $(\Ho,\Hx)$ is PE. However, the weaker condition that the pair $(\Ho,\Ho)$ is PE applies in some special situations, for example, when the trajectory $t\mapsto x(t)$ is contained in a limit cycle. This special situation is discussed in \cite{jiaACC}. As another example, When the index set $\Omega$ is finite-dimensional, which we express as $\Omega_n$, the functions in $\Hn$ can be expressed as $f=\sum\alpha_i\phi_i$. It is not difficult to see that the pair $(\Hn,\Hn)$ being PE is equivalent to the PE condition in $\mathbb{R}^n$ as in Definition \ref{def:PE_rn}.

\section{Convergence of RKHS Embedding Method}
\label{sec:convergence}
This section develops the primary convergence result in this paper that holds when the closed RKH subspace $\Ho$ is partially PE. The proof of Theorem \ref{th:pe_uas} is in the spirit of the proofs in \cite{narendra1977b} for finite-dimensional systems or the ones in \cite{demetriou1997} for infinite-dimensional systems. Here the proof is adapted to exploit the uniform bound $\|\calE^*_x\|_{op}\leq \bar{k}$ in RKHS.
\begin{theorem}
\label{th:pe_uas}
Suppose the unknown function $f\in\Hx$, where $\Hx$ is an RKHS uniformly embedded in $\calC(X)$. If the trajectory $t\mapsto x(t)$ partially persistently excites the subset $\Omega$ and the corresponding RKH subspace $\Ho$ in the sense of Definition \ref{def:PE1} (PE.1), then we have
\begin{equation*}
    \lim_{t\rightarrow \infty} \|\til{x}\| = 0, \quad \text{and} \quad \lim_{t\rightarrow \infty} \|\prj \til{f}\|_{\Hx} = 0,
\end{equation*}
for any initial condition $(\tilde{x}_0,\tilde{f}_0)$. Moreover, we have the following pointwise convergence
\begin{equation*}
    \lim_{t\rightarrow\infty}|\tilde{f}(t,x)| = 0,\quad \text{for all }x\in\Omega.
\end{equation*}
\end{theorem}

\begin{proof}
Without loss of generality, we assume that $\mu=1$ in the error Equation \eqref{eq:err_rkhs0}. Consider the candidate Lyapunov function 
\begin{equation*}
V(t)=(\til{x},P\til{x})+\dotp{\til{f},\til{f}}.    
\end{equation*}
Clearly, $V(x)$ is bounded below by zero. We take the time derivative of $V(t)$ along the trajectory of the error Equation \eqref{eq:err_rkhs0}. Applying the Lyapunov equation and using the definition of the adjoint operator, we have
\begin{align*}
    \hspace{-1em}\dot{V}(t) &= \big(\til{x},(A^TP+PA)\til{x}\big) + 2\big(B\ev{t}\til{f}, P\til{x}\big) \\
    &\qquad + 2\dotp{\til{f},-(B\ev{t})^*Px}, \\
    &= -\big(\til{x},Q\til{x}\big) + 2\big(B\ev{t}\til{f}, P\til{x}\big) \\
    &\qquad - 2\big(B\ev{t}\til{f},Px\big), \\
    &= -\big(\til{x},Q\til{x}\big).
\end{align*}
Note $Q$ is positive definite, so $\dot{V}(t)\leq 0$. This implies $V(t)$ is nonincreasing. By the Lyapunov stability theorem we conclude the error Equation \eqref{eq:err_rkhs0} is stable at the origin. 

When we integrate $\dot{V}(t)$ and take the limit for both sides as $t\rightarrow\infty$, we have
\begin{equation*}
    \int_{t_0}^{\infty} \big(\tilde{x}(\tau),Q\tilde{x}(\tau)\big) d\tau = V(t_0)-\lim_{t\rightarrow\infty}V(t) \leq V(t_0).
\end{equation*}
If the integrand $\big(\tilde{x}(t),Q\tilde{x}(t)\big)$ is uniformly continuous with respect to time $t$, then by the extension of Barbalat's lemma to Banach spaces \cite{farkas2016barbalat}, it follows that $ \big(\tilde{x}(t),Q\tilde{x}(t)\big)\rightarrow0$ as $t\rightarrow\infty$, which further implies 
\begin{equation}
\label{eq:sublemma1}
    \lim_{t\rightarrow \infty} \|\tilde{x}(t)\|=0.
\end{equation}
To show $\big(\tilde{x}(t),Q\tilde{x}(t)\big)$ is uniformly continuous, it suffices to show it is Lipschitz continuous, or equivalently that $\|\dtil{x}\|$ is uniformly bounded. From Equation \eqref{eq:err_rkhs0}, we have
\begin{equation}
\label{eq:norm_err_rkhs0}
    \|\dtil{x}\|\leq \|A\|_{op}\|\tilde{x}(t)\| + \|B\mathcal{E}_{x(t)}\|_{op}\|\tilde{f}(t)\|_{\Hx}.
\end{equation}
By the assumption of $\Hx\hookrightarrow\mathcal{C}(X)$ stated in Section \ref{sec:convergence}, the operator $\calE_x$ is uniformly bounded, i.e.  $\|\mathcal{E}_x\|_{op}\leq \bar{k}$. Due to the finiteness of $V(t)$ at each $t\geq t_0$, we have $\tilde{x}\in \mathcal{L}^{\infty}([t_0,\infty),\Rd)$, and $\tilde{f}\in \mathcal{L}^{\infty}([t_0,\infty),\Hx)$. Thus $\|\til{x}\|$ and $\|\til{f}\|_{\Hx}$ are both uniformly bounded. Thus in Equation \eqref{eq:norm_err_rkhs0}, every term is uniformly bounded, and Barbalat's lemma applies.

Now we prove the convergence of function estimate. By Equation \eqref{eq:sublemma1}, for any $\epsilon>0$, there exists $T$ such that for all $t\geq T$, $\|\til{x}\|<\epsilon$. Consider the function error $\tilde{f}(T)$ at time $T$. By the PE.1 condition, there exists $s\in [T,T+\Delta]$ such that 
\begin{equation*}
    \left| \int_s^{s+\delta} \ev{\tau} \tilde{f}(T) d\tau \right| \geq \gamma\|\prj\tilde{f}(T)\|_{\Hx}.
\end{equation*}
At time $(s+\delta)$, the state error $\|\tilde{x}(s+\delta)\|$ is bounded below as 
\begin{align}
\label{eq:term0}
    &\hspace{-2em}\|\tilde{x}(s+\delta)\| = \left\|\tilde{x}(s) + \int_{s}^{s+\delta} A\tilde{x}(\tau) + B\ev{\tau}\tilde{f}(\tau) d\tau \right\| \nonumber \\
    &\hspace{-2em}\geq \underbrace{\left\|\int_{s}^{s+\delta} B\ev{\tau}\tilde{f}(T) d\tau\right\|}_{\text{term 1}} - \underbrace{\left\|\tilde{x}(s) + \int_{s}^{s+\delta} A\tilde{x}(\tau) d\tau\right\|}_{\text{term 2}} \nonumber \\ 
    &\hspace{4em}- \underbrace{\left\|\int_{s}^{s+\delta} B\ev{\tau}(\tilde{f}(\tau)-\tilde{f}(T)) d\tau\right\|}_{\text{term 3}}.
\end{align}
In term 1 we use the PE.1 condition for $\tilde{f}(T)$, which results in the lower bound
\begin{align}
\label{eq:term1}
    \text{term 1} &= \|B\| \left|\int_{s}^{s+\delta} \ev{\tau}\tilde{f}(T) d\tau\right|,  \nonumber \\
    &\geq \gamma_1\|B\|\|\prj \tilde{f}(T)\|_{\Hx}.
\end{align}
In term 2, noting that $\|\tilde{x}(t)\|<\epsilon$ for all $t\geq T$, we have
\begin{align}
\label{eq:term2}
    \text{term 2} & \leq \|\tilde{x}(s)\| + \int_s^{s+\delta} \|A\|_{op} \|\tilde{x}(\tau)\| d\tau, \nonumber \\
    &\leq \epsilon + \epsilon\delta\|A\|_{op}.
\end{align}
For term 3, we first derive the bound on $\tilde{f}(\tau)-\tilde{f}(T)$ for $\tau\in[T,T+\Delta]$. Again, the uniform bound on the evaluation operator $\|\mathcal{E}_x\|_{op}\leq\bar{k}$ is used here.
\begin{align*}
    \|\tilde{f}(\tau)-\tilde{f}(T)\|_{\Hx} &= \left\| \int_T^\tau (B\ev{\xi})^*P\tilde{x}(\xi) d\xi \right\|_{\Hx}, \\
    &\hspace{-4em}\leq \int_T^\tau \|B^*\|_{op} \|\ev{\xi}^*\|_{op} \|P\|_{op} \|\tilde{x}(\xi)\| d\xi, \\
    &\hspace{-4em}\leq \epsilon\bar{k}(\tau-T)\|B\| \|P\|_{op}.
\end{align*}
Let $c_1=\|B\| \|P\|_{op}$, then $\|\tilde{f}(\tau)-\tilde{f}(T)\|_{\Hx}\leq c_1\bar{k}\epsilon(\tau-T)$. In term 3, note that $T\leq s\leq T+\Delta$. This means
\begin{align}
\label{eq:term3}
    \text{term 3} &\leq \int_{s}^{s+\delta} \|B\| \|\ev{\tau}\|_{op} \|\tilde{f}(\tau)-\tilde{f}(T))\|_{\Hx} d\tau, \nonumber \\
    &\leq \bar{k}\|B\| \int_{s}^{s+\delta} \epsilon c_1\bar{k}(\tau-T) d\tau, \nonumber \\
    &\leq \epsilon c_1 \bar{k}^2 \|B\|\left( \delta^2/2 + (s-T)\delta\right), \nonumber \\
    &\leq \epsilon c_1 \bar{k}^2 \|B\|\left( \delta^2/2 + \Delta\delta\right).
\end{align}
Let $c_2 = \|B\|\bar{k}^2 \left(\delta^2/2 + \Delta\delta\right)$. Then we have term 3 $\leq c_1c_2\epsilon$. Substituting Equation \eqref{eq:term1}-\eqref{eq:term3} into Equation \eqref{eq:term0} gives a lower bound for $\tilde{x}(s+\delta)$,
\begin{equation}
\label{eq:lowerbnd}
\|\tilde{x}(s+\delta)\| \geq \gamma\|B\|\|\prj\tilde{f}(T)\|_{\Hx} - (1+\|A\|_{op}\delta)\epsilon - c_1c_2\epsilon.
\end{equation}
On the other hand, we also know 
\begin{equation}
\label{eq:upperbnd}
    \|\tilde{x}(s+\delta)\|< \epsilon.
\end{equation}
When we combine the Equation \eqref{eq:lowerbnd}-\eqref{eq:upperbnd}., we get the following Equation \eqref{eq:fbound} that gives an upper bound on $\|\prj\tilde{f}(T)\|_{\Hx}$,
\begin{equation}
\|\prj\tilde{f}(T)\|_{\Hx}<\frac{\epsilon}{\gamma_1\|B\|} \big(2+\|A\|_{op}\delta + c_1c_2\big). \label{eq:fbound}
\end{equation}
Now we have shown that $\|\prj\tilde{f}(T)\|_{\Hx}$ is of $\mathcal{O}(\epsilon)$ for some $T$ that depends on $\epsilon$. It follows that $\|\prj \tilde{f}(T')\|$ is $\mathcal{O}(\epsilon)$ for all $T'\geq T$. To see why this is so, choose any $T'>T$. It is still true that $\|\tilde{x}(\tau)\|<\epsilon$ for all $\tau\geq T'>T$. We can repeat all of the steps above for $\tau\geq T'$ to conclude that $\|\prj\tilde{f}(T')\|_{\Hx}$ is $\mathcal{O}(\epsilon)$. From this we eventually conclude that
\begin{equation*}
    \lim_{t\rightarrow \infty}|\prj\tilde{f}(t)|_{\Hx}=0.
\end{equation*}
Finally we examine the evaluation of $\tilde{f}(t)$ in the PE set $\Omega$. In the view of Corollary \ref{cor:projection}, we have for all $y\in\Omega$, 
\begin{align*}
    \tilde{f}(t,y) &= \dotp{\prj\tilde{f}(t),\knl_y} + \dotp{(I-\prj)\tilde{f}(t),\knl_y}\\
    &= \dotp{\prj\tilde{f}(t),\knl_y}.
\end{align*}
Again, due to the uniform embedding $\Hx\hookrightarrow\calC(X)$, we know $\|\knl_x\|_{\Hx}\leq \bar{k}$ uniformly for all $x\in X$. Therefore, we have for all $y\in\Omega$,
\begin{equation*}
    |\tilde{f}(t,y)| \leq \bar{k}\|\prj\tilde{f}(t)\|_{\Hx},
\end{equation*}
which leads to the following pointwise convergence
\begin{equation*}
    \lim_{t\rightarrow\infty} |\tilde{f}(t,x)| = 0\quad\text{for all }x\in\Omega.
\end{equation*}
\end{proof}
\begin{remark}
Specifically, the theorem above states that when the pair $(\Ho,\Hx)$ is partially PE, we are guaranteed the convergence of RKHS embedding method in the most general situation, that is, when the unknown function $f\in\Hx$. As the price of generality, the partial PE condition tests all the functions $g\in\Hx$ with $\|\prj g\|_{\Hx}=1$ which constitute an unbounded set in $\Hx$. The PE condition can be relaxed by testing the Equation \eqref{eq:PE1} or \eqref{eq:PE2} for the functions in the subspace (i.e. $g\in\Ho$) rather than the whole RKHS $\Hx$. In that case, the theorem above is still valid with the change of assumption, in which it must holds true that $f\in\Ho\subseteq\Hx$. The conclusion can be proven following the similar steps as above. See \cite{jiaACC} for the discussions for this restricted case. 

\end{remark}

%
%
\section{Equivalence of Partial PE Conditions}
\label{sec:equiv_PE}
To motivate this section, we first recall a recent result concerning the PE condition defined for RKHS and a positive limit set $\omega^+(x_0)$.
\begin{theorem}[Theorem 2, \cite{kgp2019}]
\label{th:limit_set}
Let $\Hx$ be an RKHS of functions over $X$, and suppose $\Hx$ includes a rich family of bump functions. If the PE condition in Definition  \ref{def:PE2} (PE.2) holds for a subset $\Omega\subseteq X$, then we have $\Omega\subseteq  \omega^+(x_0).$
\end{theorem}


\noindent The theorem above states a necessary condition for a subspace $\Ho$ (or a subset $\Omega$) to be persistently excited in the sense of PE.2, for a specific type of reproducing kernel. It leads to a criterion for verifying the partial PE in practice. Given the RKHS $\Hx$ that contains a rich family of bump functions, if the subset $\Omega\subseteq X$ is not a subset of the positive limit set $\omega^+(x_0)$ of the original system, then the corresponding RKH subspace $\Ho$ cannot be persistently excited. Notice that this theorem is derived for PE.2, while in the last section the convergence of RKHS embedding method is guaranteed under PE.1. The following two theorems are devoted to exploring the connection between the two PE conditions, so as to bridge the convergence results and the practical criterion embodied in Theorem \ref{th:limit_set}.

\begin{theorem}
\label{th:PE1_implies_PE2}
The PE condition in Definition \ref{def:PE1} (PE.1) implies the one in Definition \ref{def:PE2} (PE.2).
\end{theorem}
\begin{proof}
By Definition \ref{def:PE1} (PE.1), there exist constants $T_0, \gamma_1, \delta$ and $\Delta$, such that for each $t\geq T_0$ and all $g\in \Hx$
, there is an $s\in [t,t+\Delta]$ such that 
\begin{equation*}
    \left| \int_s^{s+\delta} \ev{\tau} g d\tau \right| \geq \gamma_1 \|\prj g\|_{\Hx} >0.
\end{equation*}
Without loss of generality, we assume the interval $[s,s+\delta]$ is contained in $[t,t+\Delta]$. The integral in Equation \eqref{eq:PE2} equals
\begin{align*}
    \int_{t}^{t+\Delta}&\dotp{\ev{\tau}^*\ev{\tau}g,g} d\tau \\
    &= \int_{t}^{t+\Delta}\left(\ev{\tau}g,\ev{\tau}g \right)d\tau
    \geq  \int_{s}^{s+\delta} \big((g\circ x)(\tau)\big)^2 d\tau.
\end{align*}
We know the trajectory $t\mapsto x(t)$ is continuous. Due to the uniform embedding $\Hx \hookrightarrow\mathcal{C}(X)$, the function $g\in\Hx$ is also continuous, and so is the composition $g\circ x$. On the compact interval $[t,t+\Delta]$ the function $g\circ x$ is bounded. Thus it is also square-integrable, that is, $g\circ x\in\mathcal{L}^2([t,t+\Delta],\mathbb{R})$. By the Cauchy-Schwarz inequality,
\begin{equation*}
\begin{aligned}
&\int_{s}^{s+\delta} 1 d\tau \int_{s}^{s+\delta} \big((g\circ x)(\tau)\big)^2 d\tau
\geq \left( \int_{s}^{s+\delta} \left|\ev{\tau}g\right| d\tau \right)^2 \\ &\hspace{4em}\geq \left| \int_{s}^{s+\delta} \ev{\tau}g d\tau \right|^2 \geq \gamma_1^2\|\prj g\|_{\Hx}^2,
\end{aligned}
\end{equation*}
which then implies 
\begin{equation*}
    \int_{t}^{t+\Delta}\dotp{\ev{\tau}^*\ev{\tau}g,g} d\tau \geq \frac{\gamma_1^2}{\delta}\|\prj g\|_{\Hx}^2.
\end{equation*}
Let $\gamma_2=\gamma_1^2/\delta$. Then the PE.2 condition is satisfied with the constants $\Delta,T_0$, and $\gamma_2$.
\end{proof}
The implication above is straightforward, but the derivation of the converse needs more restrictions on the RKHS $\Hx$. Both the PE conditions require the integral of $g\circ x$ to be non-zero in a certain sense. From PE.2 and Equation \eqref{eq:PE2}, it can be deduced that $g\circ x$ has non-zero evaluation at some point $\xi\in[t,t+\Delta]$. Although $g\circ x$ is continuous, if it oscillates fast enough, the integral over a fixed-length interval $[s,s+\delta]$ is not guaranteed to be non-zero. Hence the crux of showing equivalence relies on restricting the possible fast oscillations of functions in $\Hx$. A similar pathology has been studied in \cite{narendra1977b} for the PE condition in Euclidean space.

Recall that a family of uniformly continuous functions $S$ is called uniformly equicontinuous if for each $\epsilon>0$, there is a $\delta$ depending solely on $\epsilon$ such that $\|x-y\|<\delta$ implies $|g(x)-g(y)|<\epsilon$ for all $g\in S$ and $x,y\in\Omega$. The definition of equicontinuity will be used to control the fast oscillation. Here we define the composition operator $\mathbb{U}$ associated to a trajectory $t\mapsto x(t)$ by $\mathbb{U}h:=h\circ x$ for any function $h$. Then in the following theorem, we give a sufficient condition that guarantees the equivalence of the two PE conditions.
\begin{theorem}
\label{th:pe2_to_pe1}
Let $\mathbb{U}$ be the composition operator associated with the original system in Equation \eqref{eq:orig_dyn_clas}. Define $ S_\Omega$ to be a collection of functions in $\Hx$ by
\begin{equation*}
     S_\Omega:= \{h\in\Hx: \|\prj h\|_{\Hx} = 1\}.
\end{equation*}
If the family of composite functions
\begin{equation*}
    \mathbb{U}( S_\Omega):=\{h\circ x: h\in\Hx,\; \|\prj h\|_{\Hx} = 1\}
\end{equation*}
is uniformly equicontinuous, then the PE condition in Definition \ref{def:PE2} (PE.2) implies the one in Definition \ref{def:PE1} (PE.1).
\end{theorem}

\begin{proof}
To derive PE.1, we need to  find in Definition \ref{def:PE1} the  constants $\delta,\Delta,T_0,\gamma_1>0$  and show that for each $t\geq T_0$ and $g\in \Hx$ that  there exists $s\in[t,t+\Delta]$  such that
\begin{equation}
\label{eq:th5_pe1}
    \left|\int_s^{s+\delta} \calE_{x(\tau)}g d\tau \right|\geq \gamma_1\|\prj g\|_{\Hx}.
\end{equation}
In this proof, we first derive such constants  for all the functions $h\in S_\Omega\subset \Hx$. Then we show these same constants also apply to all $g\in\Hx$.  

To begin suppose that  $h\in S_\Omega$, that is $h\in\Hx$ with $\|\prj h\|_{\Hx} = 1$. Suppose the PE.2 condition holds. Then there exist constants $T_0,\Delta$, and $\gamma_2$, such that for each $t\geq T_0$ and $h\in  S_\Omega$ the following inequality is true,
\begin{align*}
    &\hspace{-1em}\int_{t}^{t+\Delta}\dotp{\calE_{x(\tau)}^*\ev{\tau}h,h} d\tau \\&\hspace{1em}= \int_{t}^{t+\Delta}h\big(x(\tau)\big)^2 d\tau
    \geq \gamma_2\|\prj h\|_{\Hx}^2 = \gamma_2.
\end{align*}
For a fixed time $t$, by the mean value theorem there exists $s_h\in[t,t+\Delta]$ such that
\begin{equation*}
    \int_{t}^{t+\Delta}h\big(x(\tau)\big)^2 d\tau = h\big( x(s_h)\big)^2\Delta, 
\end{equation*} 
from which we conclude that 
\begin{equation}
    |h(x(s_h))|\geq \sqrt{\gamma_2/\Delta}.
\end{equation}
Let $\epsilon=\frac{1}{2}\sqrt{\gamma_2/\Delta}$. Since the family of composite functions $\mathbb{U}(S_\Omega)$ is uniformly equicontinuous, there exists a $\delta:=\delta(\epsilon) >0$ such that if $|\tau-\eta|\leq \delta$, then 
\begin{equation}
\label{eq:th5_ineq}
    |h(x(\tau))-h(x(\eta))|<\epsilon=\frac{1}{2}\sqrt{\gamma_2/\Delta}
\end{equation}
for any $h\circ x\in \mathbb{U}(S_\Omega)$. We can apply this inequality for the choice $\eta:=s_h$. 
By the triangle inequality, we know for all $\tau\in [s_h,s_h+\delta]$,
\begin{equation*}
    |h(x(\tau))| \geq \frac{1}{2}\sqrt{\gamma_2/\Delta}.
\end{equation*}
On the other hand, we also know from Equation \eqref{eq:th5_ineq} that in the interval $[s_h,s_h+\delta]$, the sign of $h\circ x$ does not change.
Thus it follows that
\begin{equation*}
    \left| \int_{s_h}^{s_h+\delta} \ev{\tau} h d\tau \right| = \int_{s_h}^{s_h+\delta} \left|\ev{\tau} h\right| d\tau \geq \frac{1}{2}\delta\sqrt{\gamma_2/\Delta}.
\end{equation*}
In this way, we have found constants $\delta$, $\Delta$, $T_0$, and $\gamma_1 = \frac{\delta}{2}\sqrt{\gamma_2/\Delta}$ such that for each $t\geq T_0$ and $h\in S_\Omega$ there is an $s:=s_h$  that satisfies the inequality in  Equation \eqref{eq:th5_pe1}. 

Now we prove these constants also work for all $g\in\Hx$. First of all, the Equation \eqref{eq:th5_pe1} is clearly true for all $g\in\Hx$ with $\|\prj g\|_{\Hx} = 0$. For those $g\in\Hx$ with $\|\prj g\|_{\Hx}\neq 0$, notice that $g/\|\prj g\|_{Hx}$ is an element of $S_\Omega$. By the derivations above, we know there exists an $s\in[t,t+\Delta]$ such that
\begin{equation*}
    \left|\int_s^{s+\delta} \calE_{x(\tau)} \frac{g}{\|\prj g\|_{\Hx}}\right| \geq \gamma_1,
\end{equation*}
which is equivalent to
\begin{equation*}
    \left|\int_s^{s+\delta} \calE_{x(\tau)} g\right| \geq \gamma_1\|\prj g\|_{\Hx}.
\end{equation*}
Hence,  the PE.1 condition is satisfied.
\end{proof}


\subsection{A Sufficient Condition for Uniform Equicontinuity}
The analysis in the last section shows that the uniform equicontinuity of the family $\mathbb{U}(S_\Omega)$ is enough to guarantee that the definitions of partial  persistency in PE.1 and PE.2 are equivalent.  Still,  this is a rather abstract condition, and it is not immmediately clear when such uniform equicontinuity might arise. 
In this section we give a simple analysis that guarantees the uniform equicontinuity of $\mathbb{U}(S_\Omega)$ in some cases. 

We suppose that we have the following situation, which is actually quite common in applications to the estimation of uncertain ODEs, see \cite{bmpkf2017,bmpkf2017C}. We suppose that that the original system trajectory $t\mapsto x(t)$ is the state of an evolution problem formulated in $X:=\RR^d$, but it exhibits at lot of structure: it is concentrated on a set $\Omega \subset X$.  Specifically, we suppose that the orbit $\Gamma^+(x_0):=\bigcup_{t\geq 0} x(t)$ is known to be contained in the subset $\Omega\subset X$. Most commonly $\Omega$ is bounded, often compact.  We begin our analysis of the uniform equicontinuity of $\mathbb{U}(S_\Omega)$ with a straightforward observation: the uniform continuity of the original trajectory and equicontinuity  of the restrictions of $S_\Omega$ imply the uniform equicontinuity of $\mathbb{U}(S_\Omega)$.  
\begin{theorem}
Let the trajectory $t\mapsto x(t)$ of the original system be uniformly continuous, choose $\Omega:=\Gamma^+(x_0)$,  and suppose that the family of restrictions 
$$
\restr(S_\Omega):=\{ g: \Omega \rightarrow \RR \ | \ g= \restr f, f\in S_\Omega \}
$$
is uniformly equicontinuous. Then $\mathbb{U}(S_\Omega)$ is uniformly equicontinuous. 
\end{theorem}
\begin{proof}
The proof follows from a simple application of the  definitions of uniform equicontinuity and uniform continuity.
Let $\epsilon>0$ be fixed.  Since $\restr(S_\Omega)$ is uniformly equicontinuous, we know that there is a $\gamma:=\gamma(\epsilon)$ such that \begin{equation}
     y,z\in \Omega \text{ and } \|y-z\|\leq \gamma \text{  implies } |g(y)-g(z)|\leq \epsilon
     \label{eq:z1}
\end{equation}
for all $g\in \restr(S_\Omega)$. But the uniform continuity of $t\mapsto x(t)$ means that for any $\gamma>0$ there is a $\delta=\delta(\gamma)$ such that 
\begin{equation}
    |\tau-s|\leq \delta \text{ implies } |x(\tau)-x(s)|\leq \gamma.
    \label{eq:z2}
\end{equation}
By combining Equations \ref{eq:z1} and \ref{eq:z2}, we conclude that the family $\mathbb{U}(S_\Omega)$ is uniformly equicontinuous. 
\end{proof}


The following lemma gives a condition that can be used as a guide for selecting reproducing kernels $\knl$ such that
the set of restrictions $\restr(S_{\Omega})$ is uniformly equicontinuous.
\begin{lemma}
Let $\Omega$ be a set endowed with the metric $d(\cdot,\cdot)$. Suppose $\knl:\Omega\times \Omega\rightarrow\mathbb{R}$ is a reproducing kernel that is expressed in terms of a  \textbf{symmetric radial basis}, that is, there exists a function $\calR:\mathbb{R}^+\rightarrow\mathbb{R}$ such that
\begin{equation*}
    \knl(x,y) \equiv \calR\big(d(x,y)\big) \quad\text{ for all } x,y\in \Omega.
\end{equation*}
If the function $\calR(\cdot)$ is $\alpha$-H\"{o}lder continuous, then the functions in $\restr(S_{\Omega})$ are uniformly equicontinuous.
\end{lemma}
\begin{proof}
Let $f$ be any function in $S_{\Omega}$. By the reproducing property, we have the following equations for any $x,y\in \Omega$,
\begin{align}
    |f(x)-f(y)| &= \left|\dotp{f,\knl_x-\knl_y}\right| \nonumber \\
    &\leq \left |\langle \prj f, \knl_x-\knl_y\rangle 
    +
    \langle (I-\prj) f, \knl_x-\knl_y\rangle 
    \right| \nonumber \\
     &\leq \|\prj f\|_{\Hx}\|\knl_x-\knl_y\|_{\Hx} \nonumber \\ \label{eq:diff_f}&=\|\knl_x-\knl_y\|_{\Hx}
\end{align}
Consider the norm of $\knl_x-\knl_y$. The reproducing property is also applied in the following equations,
\begin{align}
    \|\knl_x-\knl_y\|_{\Hx}^2 &= \dotp{\knl_x-\knl_y,\knl_x-\knl_y} \nonumber \\
    &=\knl(x,x) + \knl(y,y) - 2\knl(x,y) \label{eq:sym_knl} \\
    &= 2\calR(0) - 2\calR\big(d(x,y)\big), \label{eq:diff_k}
\end{align}
where the symmetry of $\knl$ is used in Equation \eqref{eq:sym_knl}. Suppose the function $\calR$ is $\alpha$-H\"{o}lder continuous. Then there exists a constant $C>0$ such that 
\begin{equation}
\label{eq:diff_r}
    \left|\calR(0) - \calR\big(d(x,y)\big)\right| \leq Cd(x,y)^\alpha. 
\end{equation}
We combine Equation \eqref{eq:diff_f}-\eqref{eq:diff_r}, and we get 
\begin{equation}
\label{eq:diff_bnd}
    |f(x)-f(y)|\leq \sqrt{2C}d(x,y)^{\alpha/2}.
\end{equation}
For any $\epsilon>0$ if we set $\delta<(\epsilon^2/2C)^{1/\alpha}$, then for all $y$ with $d(x,y)\leq\delta$, it holds $|f(x)-f(y)|<\epsilon$. Since all the arguments above do not rely on the choice of $f\in S_{\Omega}$ and $x\in \Omega$, it can be concluded that the functions in $\restr(S_{\Omega})$ are uniformly equicontinuous.
\end{proof}




\section{Approximation and Numerical Simulation}
\label{sec:num_sim}
\noindent Practical implementations of the estimator in Equation \eqref{eq:rkhs_estimator} requires finite-dimensional approximations. An initial discussion of approximations can be found  in \cite{bmpkf2017C,bmpkf2017}. For completeness, we review the prototypical setup of approximations in Section \ref{sec:fd}.  We subsequently  discuss in Section \ref{sec:choice} how the partial PE condition influences the choice of basis and the selection of the  kernel that defines $\Hx$. A discussion of how the newly introduced partial PE conditions impact an example problem is presented in Section \ref{sec:ex}.  

\subsection{Finite-Dimensional Approximation}
\label{sec:fd}
The finite-dimensional approximation of the function estimate $\hat{f}(t)$ is constructed in an RKH subspace $\Hn$ spanned by the basis functions $\{ \knl_{z_j} \}_{j=1}^n$. We denote the collection of kernel centers by $\Omega_n=\{z_j\}_{j=1}^n$, and the approximation of the function estimate by $\hat{f}_n(t)\in\Hn$. Ordinarily, the centers are selected along the trajectory $t\mapsto x(t)$, although this is not strictly necessary. If they are, the overall estimation process can be viewed as a data-driven method of estimation.  According to \cite{bmpkf2017C,bmpkf2017}, the governing equations of the approximation are written as
\begin{equation}
\label{eq:approx_estimator}
\begin{aligned}
\dot{\hat{x}}_n(t) &= A\hat{x}_n(t) + B\ev{t}\mathbf{P}^*_{\Omega_n}\hat{f}_n(t), \\
\dot{\hat{f}}_n(t) &= \mu \mathbf{P}_{\Omega_n}\left(B\ev{t}\right)^* P(x(t)-\hat{x}_n(t)),
\end{aligned}
\end{equation}
where $\mathbf{P}_{\Omega_n}:\Hx\rightarrow\Hn$ is the orthogonal projection operator onto the subspace $\Hn$, and its adjoint $\mathbf{P}^*_{\Omega_n}$ is interpreted as the canonical injection of $\mathcal{H}_{\Omega_n}$ into $\Hx$. Carefully note that the approximation $\hat{f}_n(t)$ in Equation \eqref{eq:approx_estimator} is \textit{not} equal to the projection of the infinite-dimensional estimate, $\mathbf{P}_{\Omega_n}\hat{f}(t)$. To see why this is true, we decompose the unknown function $f=\mathbf{P}_{\Omega_n}f + f_V$, then the error equation for the approximation is written as  
\begin{align*}
    \dot{\tilde{x}}_n(t) &= A\tilde{x}_n(t) + B\ev{t}\mathbf{P}^*_{\Omega_n}\tilde{f}_n(t) - B\ev{t}f_V, \\
    \dot{\tilde{f}}_n(t) &= \mu \mathbf{P}_{\Omega_n}\left(B\ev{t}\right)^* P\tilde{x}_n(t).
\end{align*}
If we compare the error equations above and the ones for the infinite-dimensional estimator,  we observe that the error equation with respect to $\tilde{x}_n$ is perturbed by $f_V$. See \cite{pgk2020a} for detailed analysis for the behavior of such an error system. Despite the existence of the perturbation, it is shown in \cite{bmpkf2017,bmpkf2017C} that as $n\rightarrow\infty$, that is, when the kernel centers asymptotically are dense the domain of interest $\Omega$, the solutions of the approximation equations converge to the ones of the infinite-dimensional equation on a compact time interval $[0,T]$. That is, as $n\rightarrow\infty$, we have
\begin{equation}
\label{eq:converge_approx}
\hat{x}_n(t)\rightarrow \hat{x}(t) \text{ in }\Rd,\quad\text{and}\quad \hat{f}_n(t)\rightarrow \hat{f}(t) \text{ in } \Ho.
\end{equation}

At each time $t$, the approximation $\hat{f}_n(t)$ can be written as the following linear-in-parameter expression
\begin{equation}
\label{eq:fnhat}
\hat{f}_n(t,\cdot)= \sum_{j=1}^n\hat{\alpha}_j(t)\knl_{z_j}(\cdot) = \hat{\bm{\alpha}}^T(t)\knl_{\mathbf{Z}}(\cdot),
\end{equation}
where $\hat{\bm{\alpha}}(t)=[\hat{\alpha}_{1}(t),...,\hat{\alpha}_{n}(t)]^T$ is the vector of time-varying coefficients, and $\knl_{\mathbf{Z}} = [\knl_{z_1},...,\knl_{z_n}]^T$ is the vector of basis functions. Using the reproducing property, the equations of finite-dimensional approximation in Equation \eqref{eq:approx_estimator} can be transformed to those with respect to the coefficients $\hat{\bm{\alpha}}(t)$ , which take the form
\begin{align*}
    \dot{\hat{x}}_n(t) &= A\hat{x}_n(t) + B\hat{\bm{\alpha}}^T(t)\knl_{\mathbf{Z}}(x(t)), \\
    \dot{\hat{\bm{\alpha}}}(t) &= \mu\mathbb{K}^{-1}\knl_{\mathbf{Z}}(x(t))B^TP(x(t)-\hat{x}_n(t)),
\end{align*}
where $\mathbb{K}=[\knl(z_i,z_j)]_{1\leq i,j\leq n}$ denotes the Grammian matrix associated with the centers $\{z_j\}_{j=1}^n$. 

If we denote the explicit form of the projection by $\mathbf{P}_{\Omega_n}f=\bm{\alpha}^*\knl_Z$, a natural question arises whether the coefficients $\bm{\alpha}(t)$ converge to the ``ideal'' values $\bm{\alpha}^*$. The answer is unfortunately no, because of the perturbation $f_V$. However, in \cite{pgk2020a} it is shown that $\bm{\alpha}(t)$ eventually converges to a neighborhood of $\bm{\alpha}^*$ as $t\rightarrow\infty$. This analysis is entirely analogous to the study of ultimate boundedness in conventional adaptive estimation in Euclidean spaces. \cite{PoFar} 


\subsection{Choice of Kernel Basis Functions}
\label{sec:choice}
The classical  PE conditions for adaptive estimation in Euclidean space are notoriously difficult to establish, or to use in constructing pragmatic estimation schemes. As we explain in this section, however, the newly introduced partial PE conditions can be important tools in guiding or constructing pragmatic estimation methods.  
Theorem \ref{th:limit_set} shows that, with a judicious choice of the kernel $\knl$, a PE index set $\Omega$ is necessarily contained in the positive limit set $\omega^+(x_0)$ of the original system. This suggests that a reasonable choice of kernel centers $\Omega_n$ should constitute a good sampling of the positive limit set $\omega^+(x_0)$ of the orbit, denoted by $\Gamma^+(x_0)=\bigcup_{\tau\geq0}x(\tau)$. In this way, we seek estimates that converge in the indexing set $\Omega_n\subseteq \Omega \equiv \omega^+(x_0)$. Of course, in practice,  the explicit form of the positive limit set is not always known. In some cases it can be estimated using Lyapunov function techniques. See \cite{pgk2020b} for a few  specific algorithms for estimating the unknown positive limit set from measurement of the states.  Thus, in contrast to our understanding of the PE in the Euclidean spaces, the new PE conditions establish  that approximation and estimation of the unknown function is  inextricably connected to identifying regions $\Omega$ that are visited repeatedly by the unknown trajectory. In this sense, the analysis of this paper makes rigorous the intuition expressed in papers like \cite{wang2006,wu2019,farrell99wavelet} where qualitative relationships between the unknown trajectory and the fidelity of approximation are described. 

To be specific, we consider a  typical case. We assume that the trajectory $t\mapsto x(t)$ evolves on an unknown smooth, compact, connected $m$-dimensional regularly embedded submanifold $M\equiv\omega^+(x_0)\subseteq \mathbb{R}^d$. While in practice the positive limit set $\Gamma^+(x_0)$ is not known, in the analysis that follows, for the purpose of studying convergence of function estimates, we assume that it is known.  We specifically employ the extrinsic methodology summarized in Corollary 4 of \cite{kgp2019} to a system that has the positive limit set $\omega^+(x_0)$.  
We choose the $\nu$-th order Sobolev-Mat\'{e}rn kernel $\knl^\nu$ on $\Rd$ to define the RKHS $\Hx$ with $X=\Rd$. This is a classical choice of kernel over $\Rd$, which induces an RKHS isometrically isomorphic to some Sobolev space $W^{\tau,2}(\RR^d)$. By the Sobolev embedding theorem, when the smoothness index $\tau>d/2$, the space $W^{\nu,2}(\RR^d)$  is embedded in the space of continuous functions. The Sobolev-Mat\'{e}rn kernel is radial basis and symmetric, and has the following closed form expression
\begin{equation}
\begin{aligned}
\label{eq:matern_kernel}
\knl^\nu(x,y) &= \calR^\nu(\|x-y\|) \\
&\hspace{-2em} = \calR^\nu(\xi) = \frac{2^{1-\nu}}{\Gamma(\nu)} \left(\frac{\sqrt{2\nu}\xi}{l}\right)^{\nu} \mathcal{B}_{\nu}\left(\frac{\sqrt{2\nu}\xi}{l}\right),
\end{aligned}
\end{equation}
where $\xi=\|x-y\|$, $l$ is the scaling factor of length, and $\mathcal{B}_{\nu}$ is the Bessel function of order $\nu$. When $\nu=p+1/2$ with a non-negative integer $p$, the kernel can be expressed as the product of a polynomial and an exponential \cite{rasmussen}. 

\subsection{Numerical Example and Simulation}
\label{sec:ex}
\begin{figure}
    \centering
    \includegraphics[width=.35\textwidth]{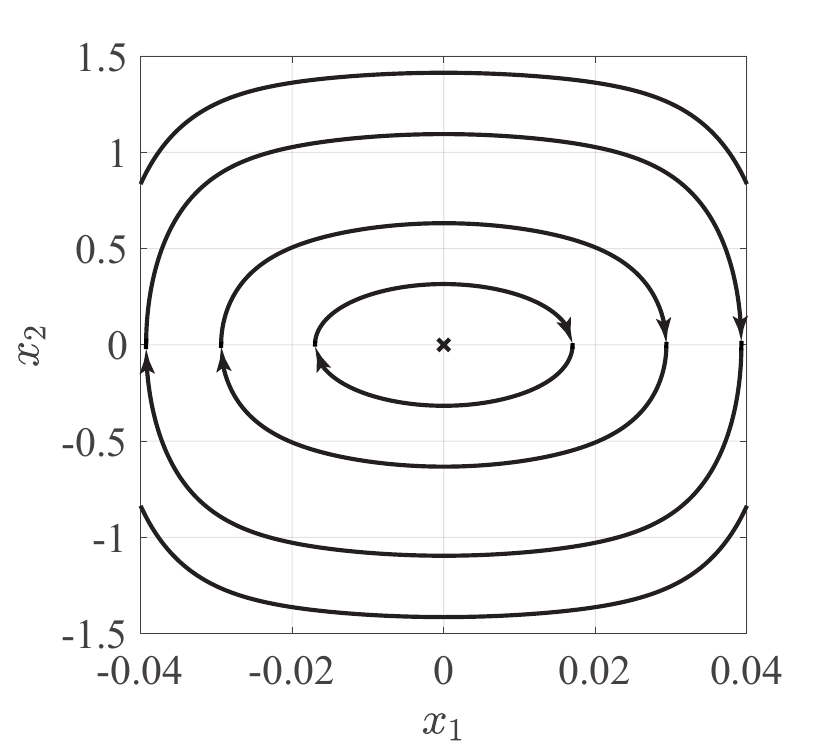}
    \caption{Phase portrait of the original system}
    \label{fig:pp}
\end{figure}

\begin{figure}
    \centering
    \includegraphics[width=.45\textwidth]{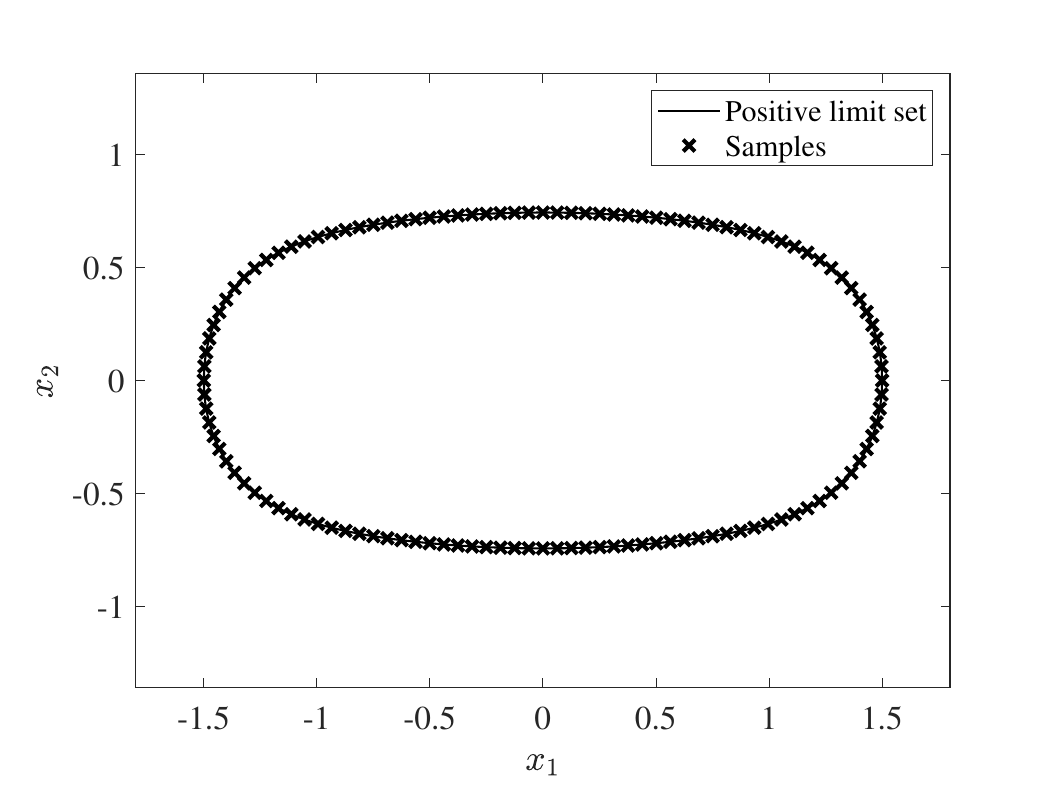}
    \caption{Samples along the positive limit set}
    \label{fig:samples}
\end{figure}

To illustrate the convergence of RKHS embedding method and the behavior of finite-dimensional approximation, an undamped, nonlinear, piezoelectric oscillator studied in \cite{sai1,sai2} is taken as an example. The governing equations of the oscillator, after employing a single bending mode approximation, have the form
\begin{align}
\label{eq:beam_dyn}
\dot{x}_1 &=  x_2, \\ 
\dot{x}_2 &=  -\frac{k}{m} x_1 - \frac{k_{n,1}}{m} x_1^3 - \frac{k_{n,2}}{m} x_1^5,
\end{align}
where $m$ is the mass, $k$ is the electromechanical stiffness, and $k_{n,1}, k_{n,2}$ are the higher order electromechanical stiffness coefficients. Figure \ref{fig:pp} shows the typical positive limit sets of this system. The limit sets take the form of limit cycles around the equilibrium at the origin, which is prototypical for such conservative electromechanical oscillators. For this example, the initial condition is selected such that the orbit $\Gamma^+(x_0)$ is precisely the limit cycle shown in Figure\ref{fig:samples}. 

In the governing equations above, we assume all the linear terms are known, and the nonlinear term 
$$f(x) = - \frac{k_{n,1}}{m} x_1^3 - \frac{k_{n,2}}{m} x_1^5$$
is unknown and to be estimated. In this case, the Sobolev-Mat\'{e}rn kernel and the associated RKHS is applied, which is uniformly embedded in the space of continuous functions \cite{wendland}. Here we set the order $\nu=3/2$ for the form expressed in Equation \eqref{eq:matern_kernel}, which yields the following kernel
\begin{equation*}
    \knl^{3/2}(x,y) = \left(1+\frac{\sqrt{3}\|x-y\|}{l}\right)\exp\left(-\frac{\sqrt{3}\|x-y\|}{l}\right).
\end{equation*}
From the conclusion of \cite{kgp2019}, a persistently excited set $\Omega$ must be contained in the positive limit set of the system $\omega^+(x_0)$, which we approximate with carefully chosen samples $\Omega_n=\{z_j\}_{j=1}^n$. The samples $\{z_j\}$ serving as the centers of basis functions $\{\knl_{z_j}\}$ are chosen along the positive limit set. 
\begin{figure}
    \centering
    \includegraphics[width=.45\textwidth]{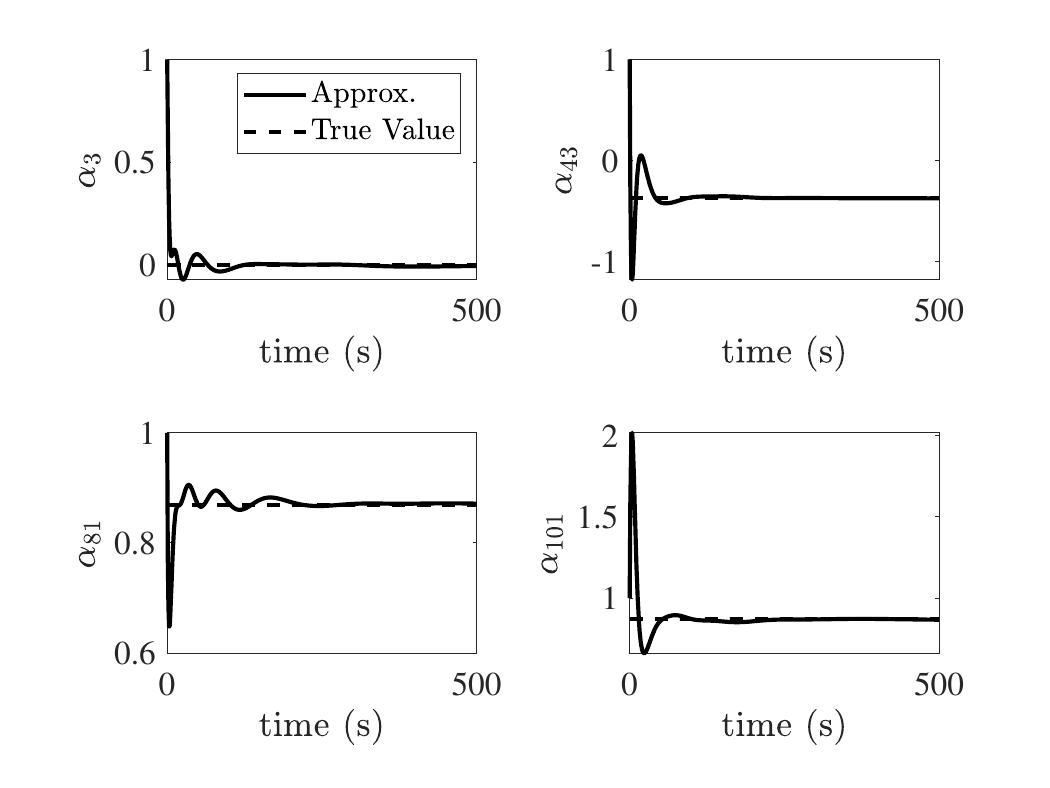}
    \caption{Convergence of parameters}
    \label{fig:param_conv}
\end{figure}

\begin{figure}
    \centering
    \includegraphics[width=.4\textwidth]{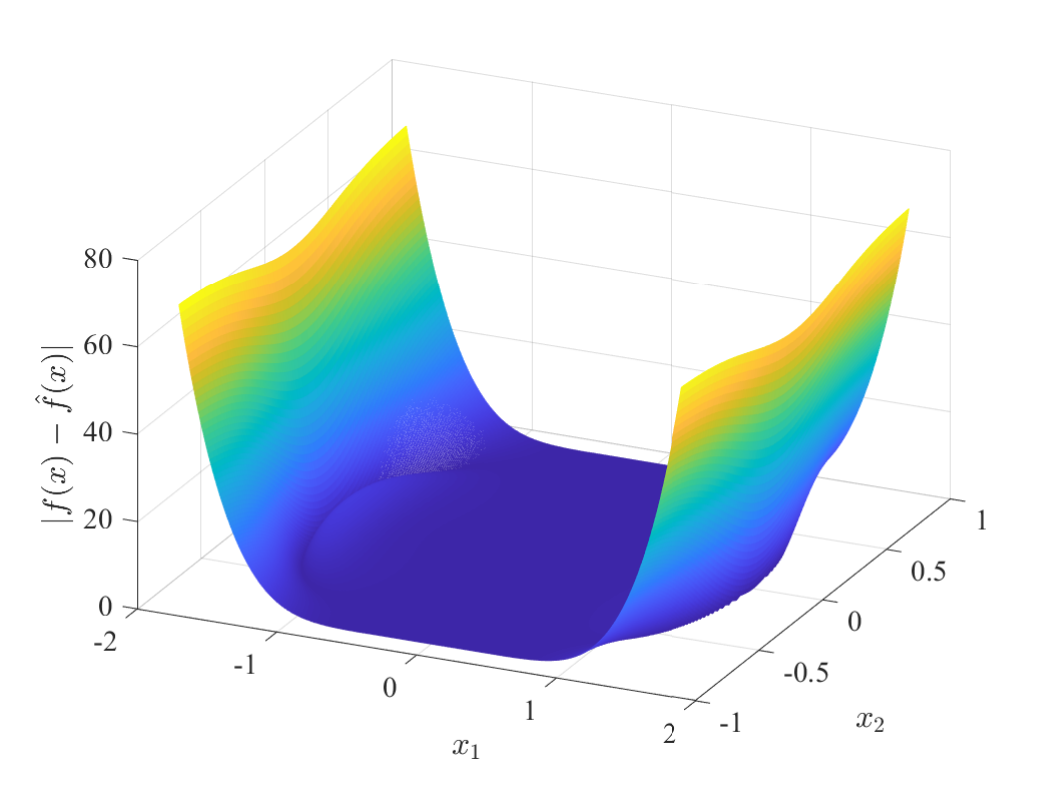}
    \caption{Error in function estimate}
    \label{fig:errsurf}
\end{figure}

\begin{figure}
    \centering
    \includegraphics[width=.4\textwidth]{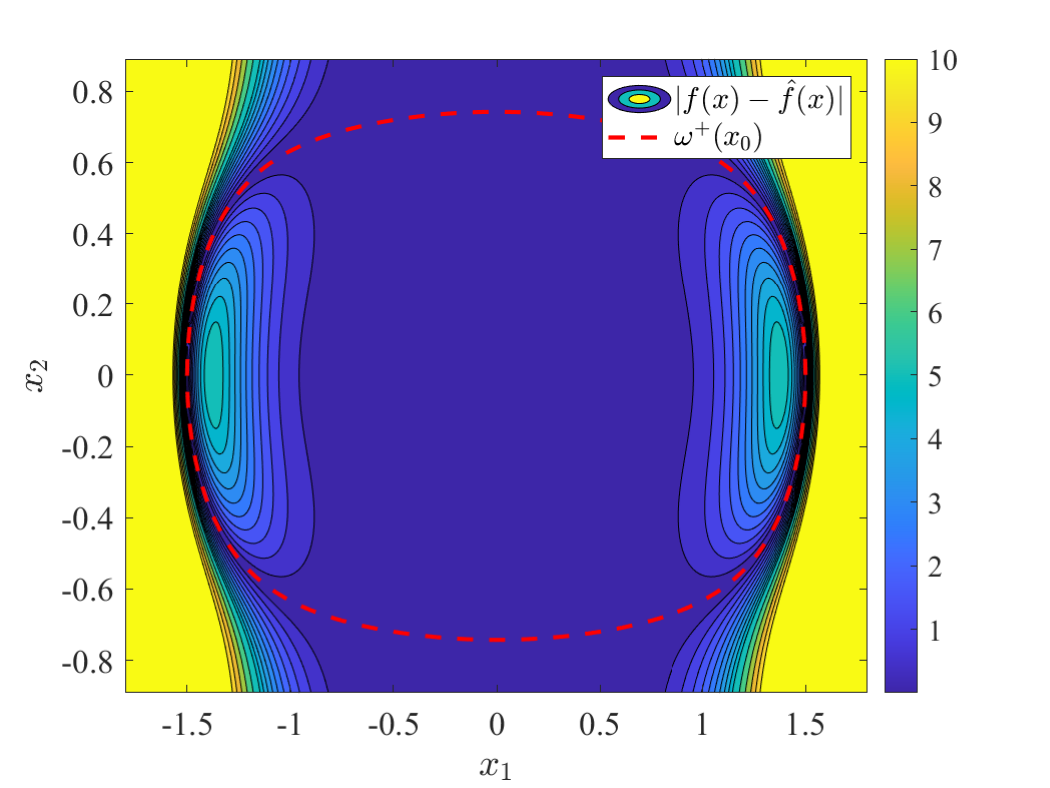}
    \caption{Error contour}
    \label{fig:errcont}
\end{figure}

When the approximation of the infinite-dimensional adaptive estimator based on the RKHS embedding technique is implemented for this problem, estimates of the unknown nonlinear function $f$ are obtained in $\Hn=\text{span}\{\knl_{z_j}\}_{j=1}^n$. Figures \ref{fig:samples}-\ref{fig:errcont} depict the results of simulation for the specific case when the number of centers $n=120$. The parameters in the state Equation \eqref{eq:beam_dyn} are $m=0.9745$ $k=6.5980$, $k_{n,1} = -1.0320$, and $k_{n,2}=3.8568$. Figure \ref{fig:param_conv} shows the time history of four typical coefficients $\hat{\alpha}_i(t)$ compared to the corresponding true value $\alpha^*_i$ calculated from projecting the true function $f$ onto the RKHS $\Hn$. There exist a slight difference between the convergent coefficient with the true value as we would expect for the case when $n$ is fixed: the coefficient estimates for fixed $n$ satisfy ultimate boundedness conditions as in the classical case. Figure \ref{fig:errsurf} shows the error between the actual function and function estimate over the state space. Qualitatively, convergence of the function estimate is expected from the theory in this paper as $n\rightarrow\infty$ over the positive limit set of the particular trajectory. Figure \ref{fig:errcont} shows the contour of the function error along with the positive limit set. Both the figures suggest that the function estimate in $\Ho$ converges to the actual function over the indexing set, which when persistently excited is a subset of the positive limit set. 

\section{Conclusions}
\noindent In this paper, new partial PE conditions are introduced for the non-parametric adaptive estimation based on RKHS embedding. The paper proves that PE.1 in Definition \ref{def:PE1} is a sufficient condition for convergence of projections defined in terms of the error equations that arise in the RKHS embedding framework. The proof relies in a fundamental way on the assumption that the RKHS is uniformly embedded in the space of continuous functions. The result provides a clear characterization of the region where the function estimate pointwise converges. Motivated by an intuitive necessary condition of PE.2 in Definition \ref{def:PE2}, this paper establishes sufficient conditions for the equivalence of PE.1 and PE.2. Condition PE.1 naturally implies PE.2, and PE.2 implies PE.1 when the family of composite functions $\mathbb{U}(S_\Omega)$ is uniformly equicontinuous. A numerical example is given to show qualitatively the convergence behavior of the RKHS embedding method. The numerical example shows that, roughly speaking, convergence is established over sets in which trajectories are eventually concentrated.

\section*{Declaration of Competing Interest}
The authors declare that they have no known competing financial interests or personal relationships that could have appeared to influence the work reported in this paper.

\section*{Acknowledgement}
Andrew J. Kurdila would like to acknowledge the support of the Army Research Office under the award \textbf{Distributed Consensus Learning for Geometric and Abstract Surfaces}, ARO Grant W911NF-13-1-0407.

\bibliographystyle{IEEEtran}
\bibliography{pe_rkhs}

\begin{thebibliography}{10}
\providecommand{\url}[1]{#1}
\csname url@samestyle\endcsname
\providecommand{\newblock}{\relax}
\providecommand{\bibinfo}[2]{#2}
\providecommand{\BIBentrySTDinterwordspacing}{\spaceskip=0pt\relax}
\providecommand{\BIBentryALTinterwordstretchfactor}{4}
\providecommand{\BIBentryALTinterwordspacing}{\spaceskip=\fontdimen2\font plus
\BIBentryALTinterwordstretchfactor\fontdimen3\font minus
  \fontdimen4\font\relax}
\providecommand{\BIBforeignlanguage}[2]{{%
\expandafter\ifx\csname l@#1\endcsname\relax
\typeout{** WARNING: IEEEtran.bst: No hyphenation pattern has been}%
\typeout{** loaded for the language `#1'. Using the pattern for}%
\typeout{** the default language instead.}%
\else
\language=\csname l@#1\endcsname
\fi
#2}}
\providecommand{\BIBdecl}{\relax}
\BIBdecl

\bibitem{PoFar}
J.~A. Farrell and M.~M. Polycarpou, \emph{Adaptive approximation based control:
  unifying neural, fuzzy and traditional adaptive approximation
  approaches}.\hskip 1em plus 0.5em minus 0.4em\relax John Wiley \& Sons, 2006,
  vol.~48.

\bibitem{slotine1991}
R.~M. Sanner and J.-J.~E. Slotine, ``Gaussian networks for direct adaptive
  control,'' in \emph{1991 American Control Conference}.\hskip 1em plus 0.5em
  minus 0.4em\relax IEEE, 1991, pp. 2153--2159.

\bibitem{slotine1992}
------, ``Stable recursive identification using radial basis function
  networks,'' in \emph{1992 American Control Conference}.\hskip 1em plus 0.5em
  minus 0.4em\relax IEEE, 1992, pp. 1829--1833.

\bibitem{sastry2011book}
S.~Sastry and M.~Bodson, \emph{Adaptive control: stability, convergence and
  robustness}.\hskip 1em plus 0.5em minus 0.4em\relax Courier Corporation,
  2011.

\bibitem{naranna}
K.~S. Narendra and A.~M. Annaswamy, \emph{Stable adaptive systems}.\hskip 1em
  plus 0.5em minus 0.4em\relax Courier Corporation, 2012.

\bibitem{IaSu}
P.~A. Ioannou and J.~Sun, \emph{Robust adaptive control}.\hskip 1em plus 0.5em
  minus 0.4em\relax Courier Corporation, 2012.

\bibitem{narendra1977b}
A.~Morgan and K.~Narendra, ``On the stability of nonautonomous differential
  equations $\dot{x}$=a+b(t)x, with skew symmetric matrix b(t),'' \emph{SIAM
  Journal on Control and Optimization}, vol.~15, no.~1, pp. 163--176, 1977.

\bibitem{anderson1977exponential}
B.~Anderson, ``Exponential stability of linear equations arising in adaptive
  identification,'' \emph{IEEE Transactions on Automatic Control}, vol.~22,
  no.~1, pp. 83--88, 1977.

\bibitem{panteley2003}
A.~Loria, E.~Panteley, D.~Popovic, and A.~R. Teel, ``Persistency of excitation
  for uniform convergence in nonlinear control systems,'' \emph{arXiv preprint
  math/0301335}, 2003.

\bibitem{annaswamy2018}
B.~M. Jenkins, A.~M. Annaswamy, E.~Lavretsky, and T.~E. Gibson, ``Convergence
  properties of adaptive systems and the definition of exponential stability,''
  \emph{SIAM journal on control and optimization}, vol.~56, no.~4, pp.
  2463--2484, 2018.

\bibitem{sastry2019}
K.~Nar and S.~S. Sastry, ``Persistency of excitation for robustness of neural
  networks,'' \emph{arXiv preprint arXiv:1911.01043}, 2019.

\bibitem{annaswamy2019}
J.~E. Gaudio, A.~M. Annaswamy, E.~Lavretsky, and M.~A. Bolender, ``Parameter
  estimation in adaptive control of time-varying systems under a range of
  excitation conditions,'' \emph{arXiv preprint arXiv:1911.03810}, 2019.

\bibitem{wang2006}
C.~Wang and D.~J. Hill, ``Learning from neural control,'' \emph{IEEE
  Transactions on Neural Networks}, vol.~17, no.~1, pp. 130--146, 2006.

\bibitem{wu2019}
W.~Wu, C.~Wang, and C.~Yuan, ``Deterministic learning from sampling data,''
  \emph{Neurocomputing}, vol. 358, pp. 456--466, 2019.

\bibitem{farrell99wavelet}
N.~Sureshbabu and J.~A. Farrell, ``Wavelet-based system identification for
  nonlinear control,'' \emph{IEEE Transactions on Automatic Control}, vol.~44,
  no.~2, pp. 412--417, 1999.

\bibitem{wendland}
H.~Wendland, \emph{Scattered data approximation}.\hskip 1em plus 0.5em minus
  0.4em\relax Cambridge university press, 2004, vol.~17.

\bibitem{rasmussen}
C.~K. Williams and C.~E. Rasmussen, \emph{Gaussian processes for machine
  learning}.\hskip 1em plus 0.5em minus 0.4em\relax MIT press Cambridge, MA,
  2006, vol.~2, no.~3.

\bibitem{aronszajn1950}
N.~Aronszajn, ``Theory of reproducing kernels,'' \emph{Transactions of the
  American mathematical society}, vol.~68, no.~3, pp. 337--404, 1950.

\bibitem{berlinet}
A.~Berlinet and C.~Thomas-Agnan, \emph{Reproducing kernel Hilbert spaces in
  probability and statistics}.\hskip 1em plus 0.5em minus 0.4em\relax Springer
  Science \& Business Media, 2011.

\bibitem{chowdhary2014}
G.~Chowdhary, H.~A. Kingravi, J.~P. How, and P.~A. Vela, ``Bayesian
  nonparametric adaptive control using gaussian processes,'' \emph{IEEE
  transactions on neural networks and learning systems}, vol.~26, no.~3, pp.
  537--550, 2014.

\bibitem{Pillonetto2010a}
G.~Pillonetto and G.~{De Nicolao}, ``{A new kernel-based approach for linear
  system identification},'' \emph{Automatica}, vol.~46, no.~1, pp. 81--93, Jan
  2010.

\bibitem{Pillonetto2011b}
G.~Pillonetto, M.~H. Quang, and A.~Chiuso, ``{A new kernel-based approach for
  nonlinearsystem identification},'' \emph{IEEE Transactions on Automatic
  Control}, vol.~56, no.~12, pp. 2825--2840, dec 2011.

\bibitem{Pillonetto2011}
G.~Pillonetto and G.~{De Nicolao}, ``{Kernel selection in linear system
  identification Part I: A Gaussian process perspective},'' in
  \emph{Proceedings of the IEEE Conference on Decision and Control}, 2011, pp.
  4318--4325.

\bibitem{Chen2018}
T.~Chen, ``{Continuous-Time DC Kernel - A Stable Generalized First-Order Spline
  Kernel},'' \emph{IEEE Transactions on Automatic Control}, vol.~63, no.~12,
  pp. 4442--4447, dec 2018.

\bibitem{Libera2019}
\BIBentryALTinterwordspacing
A.~D. Libera, R.~Carli, and G.~Pillonetto, ``{A novel Multiplicative Polynomial
  Kernel for Volterra series identification},'' \emph{arXiv preprint
  1905.07960}, 2019. [Online]. Available: \url{http://arxiv.org/abs/1905.07960}
\BIBentrySTDinterwordspacing

\bibitem{bmpkf2017C}
P.~Bobade, S.~Majumdar, S.~Pereira, A.~J. Kurdila, and J.~B. Ferris, ``Adaptive
  estimation in reproducing kernel hilbert spaces,'' in \emph{2017 American
  Control Conference (ACC)}.\hskip 1em plus 0.5em minus 0.4em\relax IEEE, 2017,
  pp. 5678--5683.

\bibitem{bmpkf2017}
------, ``Adaptive estimation for nonlinear systems using reproducing kernel
  hilbert spaces,'' \emph{Advances in Computational Mathematics}, vol.~45,
  no.~2, pp. 869--896, 2019.

\bibitem{bsdr1997}
J.~Baumeister, W.~Scondo, M.~Demetriou, and I.~Rosen, ``On-line parameter
  estimation for infinite dimensional dynamical systems,'' \emph{SIAM Journal
  of Control and Optimisation}, vol.~35, no.~2, pp. 678--713, 1997.

\bibitem{bdrr1998}
M.~B{\"o}hm, M.~Demetriou, S.~Reich, and I.~Rosen, ``Model reference adaptive
  control of distributed parameter systems,'' \emph{SIAM Journal on Control and
  Optimization}, vol.~36, no.~1, pp. 33--81, 1998.

\bibitem{d1993}
M.~A. Demetriou, ``Adaptive parameter estimation of abstract parabolic and
  hyperbolic distributed parameter systems.'' Ph.D. dissertation, University of
  Southern California, 1994.

\bibitem{dr1994}
M.~Demetriou and I.~Rosen, ``Adaptive identification of second-order
  distributed parameter systems,'' \emph{Inverse Problems}, vol.~10, no.~2, p.
  261, 1994.

\bibitem{dr1994pe}
------, ``On the persistence of excitation in the adaptive identification of
  distributed parameter systems,'' \emph{IEEE Transactions on Automatic
  Control}, vol.~39, pp. 1117--1123, 1994.

\bibitem{kl2013}
A.~Kurdila and Y.~Lei, ``Adaptive control via embedding in reproducing kernel
  hilbert spaces,'' in \emph{2013 American Control Conference}.\hskip 1em plus
  0.5em minus 0.4em\relax IEEE, 2013, pp. 3384--3389.

\bibitem{bobade2}
P.~Bobade, D.~Panagou, and A.~J. Kurdila, ``Multi-agent adaptive estimation
  with consensus in reproducing kernel hilbert spaces,'' in \emph{2019 18th
  European Control Conference (ECC)}.\hskip 1em plus 0.5em minus 0.4em\relax
  IEEE, 2019, pp. 572--577.

\bibitem{kur95}
A.~Kurdila, F.~J. Narcowich, and J.~D. Ward, ``Persistency of excitation in
  identification using radial basis function approximants,'' \emph{SIAM journal
  on control and optimization}, vol.~33, no.~2, pp. 625--642, 1995.

\bibitem{kgp2019}
A.~J. Kurdila, J.~Guo, S.~T. Paruchuri, and P.~Bobade, ``Persistence of
  excitation in reproducing kernel hilbert spaces, positive limit sets, and
  smooth manifolds,'' \emph{arXiv preprint: 1909.12274}, 2019.

\bibitem{pgk2020a}
S.~T. Paruchuri, J.~Guo, and A.~Kurdila, ``Sufficient conditions for parameter
  convergence over embedded manifolds using kernel techniques,'' \emph{arXiv
  preprint arXiv:2009.02866}, 2020.

\bibitem{Boyd1983On}
S.~Boyd and S.~Sastry, ``On parameter convergence in adaptive control,''
  \emph{Systems \& control letters}, vol.~3, no.~6, pp. 311--319, 1983.

\bibitem{jiaACC}
J.~{Guo}, S.~T. {Paruchuri}, and A.~J. {Kurdila}, ``Persistence of excitation
  in uniformly embedded reproducing kernel hilbert (rkh) spaces,'' in
  \emph{2020 American Control Conference (ACC)}, 2020, pp. 4539--4544.

\bibitem{demetriou1997}
J.~Baumeister, W.~Scondo, M.~Demetriou, and I.~Rosen, ``On-line parameter
  estimation for infinite-dimensional dynamical systems,'' \emph{SIAM Journal
  on Control and Optimization}, vol.~35, no.~2, pp. 678--713, 1997.

\bibitem{farkas2016barbalat}
B.~Farkas and S.-A. Wegner, ``Variations on barb{\u{a}}lat's lemma,'' \emph{The
  American Mathematical Monthly}, vol. 123, no.~8, pp. 825--830, 2016.

\bibitem{pgk2020b}
S.~T. Paruchuri, J.~Guo, and A.~Kurdila, ``Kernel center adaptation in the
  reproducing kernel hilbert space embedding method,'' \emph{arXiv preprint
  arXiv:2009.02867}, 2020.

\bibitem{sai1}
S.~T. Paruchuri, J.~Sterling, V.~V.~S. Malladi, A.~Kurdila, J.~Vignola, and
  P.~Tarazaga, ``Passive piezoelectric subordinate oscillator arrays,''
  \emph{Smart Materials and Structures}, vol.~28, no.~8, p. 085046, 2019.

\bibitem{sai2}
S.~T. Paruchuri, J.~Sterling, A.~Kurdila, and J.~Vignola, ``Piezoelectric
  composite subordinate oscillator arrays and frequency response shaping for
  passive vibration attenuation,'' in \emph{2017 IEEE Conference on Control
  Technology and Applications (CCTA)}.\hskip 1em plus 0.5em minus 0.4em\relax
  IEEE, 2017, pp. 702--707.

\end{thebibliography}

\section*{Appendix}
\subsection{Proof of Corollary \ref{cor:projection}}
\label{app:cor1}
\begin{proof}
Let $\Vo$ be all the functions in $\Hx$ whose restrictions over $\Omega$ are zero,
\begin{equation*}
    \Vo := \{v\in\Hx: v(x)=0 \; \forall x\in\Omega \} = \restr^{-1}\{0\}.
\end{equation*}
It is clear from its definition that  $\Vo$ is a subspace of $\Hx$. Let $\{v_n\}$ be a  sequence in $\Vo$ that converges to $v\in\Hx$. Since the norm convergence in $\Hx$ implies pointwise convergence, it follows that
\begin{equation*}
    v(x)=\lim_{n\rightarrow\infty}v_n(x) = 0,\quad \forall x\in\Omega.
\end{equation*} 
Thus we have $v\in\Vo$, and $\Vo$ is a closed subspace of $\Hx$. We define $\Ho=\Vo^\perp$ to be the orthogonal subspace of $\Vo$. Any function $f\in\Hx$ has an orthogonal decomposition 
\begin{equation*}
    f = f_\Omega + f_V,
\end{equation*}
where $f_\Omega\in\Ho$ and $f_V\in\Vo$. By the Pythagorean theorem, it holds that
\begin{equation}
\label{eq:pythagoras}
    \|f\|_{\Hx}^2 = \|f_\Omega\|_{\Hx}^2 + \|f_V\|_{\Hx}^2.
\end{equation} 
It is immediate that $\restr:\Hx\rightarrow \restr(\Hx)$ is a linear, onto mapping.  As a consequence, $\restr |_{\Ho}:\Ho \rightarrow \restr \Hx$ is linear, one-to-one,  and onto since $\Vo=\text{nullspace}(\restr)$. 
We define the extension operator $\exten:=\left (\restr|_{\Ho} \right)^{-1}$.
 It is clear that for any $f_\Omega\in\Ho$, the function $\exten\restr f_\Omega = f_\Omega$, which means the composition $\prj = \exten\restr$ is in fact a projection operator.
It is stated in Theorem \ref{th:ber6} that $\knl_\Omega:=\knl |_{\Omega \times \Omega}$ is the reproducing kernel of the space of restrictions $\restr\Hx$. By Moore-Aronszajn theorem \cite{aronszajn1950,berlinet}, the space $\restr\Hx$ is equivalent to
\begin{equation*}
    \restr\Hx = \overline{span\{\knl_{\Omega,x}: x\in\Omega \}}.
\end{equation*}
It is easy to show that $\exten\knl_{\Omega,x} = \knl_x$ for all $x\in\Omega$. Therefore, the subspace $\Ho$ can be expressed as
\begin{equation*}
    \Ho = \overline{span\{\knl_x:x\in\Omega\}}.
\end{equation*}
\end{proof}

\end{document}